\makeatletter 

\expandafter\ifx\csname amscd.sty\endcsname\relax
\expandafter\def\csname amscd.sty\endcsname{}
\else\endinput\fi
\def\filename{amscd.sty}
\def\fileversion{1.1} \def\filedate{21-JUN-1991}
\immediate\write16{%
AMS-Latex option `\filename' (\fileversion, \filedate)}
\def\Invalid@@{Invalid use of \string}
\def\Let@{\let\\\math@cr}
\def\RIfM@{\relax\protect\ifmmode}
\@ifundefined{math@cr}
  {\def\math@cr{{\ifnum0=`}\fi
   \new@ifstar{\global\@eqpen\@M\math@cr@}%
          {\global\@eqpen\interdisplaylinepenalty \math@cr@}}}
  {}
\def\math@cr@{\new@ifnextchar[\math@cr@@{\math@cr@@[\z@]}}
\def\math@cr@@[#1]{\ifnum0=`{\fi}\math@cr@@@
  \noalign{\vskip#1\relax}}
\def\restore@math@cr{\def\math@cr@@@{\cr}}
\restore@math@cr
\def\new@ifnextchar#1#2#3{%
  \let\@tempe #1\def\@tempa{#2}\def\@tempb{#3}\futurelet
    \@tempc\new@ifnch}
\def\new@ifnch{\ifx\@tempc \@tempe \let\@tempd\@tempa
             \else\let\@tempd\@tempb\fi\@tempd}
\def\new@ifstar#1#2{\new@ifnextchar *{\def\@tempa*{#1}\@tempa}{#2}}
\def\DN@{\def\next@}
\def\FN@{\futurelet\next}
\def\setboxz@h{\setbox\z@\hbox}
\def\wdz@{\wd\z@}
\def\setbox@ne{\setbox\@ne}
\def\wd@ne{\wd\@ne}
\def\rightarrowfill@#1{\m@th\setboxz@h{$#1-$}\ht\z@\z@
  $#1\copy\z@\mkern-6mu\cleaders
  \hbox{$#1\mkern-2mu\box\z@\mkern-2mu$}\hfill
  \mkern-6mu\mathord\rightarrow$}
\def\leftarrowfill@#1{\m@th\setboxz@h{$#1-$}\ht\z@\z@
  $#1\mathord\leftarrow\mkern-6mu\cleaders
  \hbox{$#1\mkern-2mu\copy\z@\mkern-2mu$}\hfill
  \mkern-6mu\box\z@$}
\def\leftrightarrowfill@#1{\m@th\setboxz@h{$#1-$}\ht\z@\z@
  $#1\mathord\leftarrow\mkern-6mu\cleaders
  \hbox{$#1\mkern-2mu\box\z@\mkern-2mu$}\hfill
  \mkern-6mu\mathord\rightarrow$}
\long\def\@leftmark#1#2{#1}
\long\def\@rightmark#1#2{#2}
\long\def\@ifempty#1{%
 \expandafter\ifx\@car#1@\@nil @\@empty
  \expandafter\@leftmark\else\expandafter\@rightmark\fi}
\long\def\@ifnotempty#1{\@ifempty{#1}{}}
\def\atdef@#1{\expandafter\def\csname\space @\string#1\endcsname}
\begingroup \catcode`\@=\active
\xdef @{\expandafter\noexpand\csname FN\string @\endcsname
  \expandafter\noexpand\csname at\string @\endcsname}
\endgroup
\def\at@{\let\next@\at@@
 \ifcat\noexpand\next a\else
 \ifcat\noexpand\next0\else
 \ifcat\noexpand\next\relax\else
 \let\next@\at@@@\fi\fi\fi\next@}
\def\at@@#1{\expandafter
  \ifx\csname\space @\string#1\endcsname\relax
    \DN@{\at@@@#1}%
  \else
    \DN@{\csname\space @\string#1\endcsname}%
  \fi\next@}%
\def\at@@@{\err@{\Invalid@@ @}{\the\athelp@}\char64\relax}
\@ifundefined{athelp@}{\csname newhelp\endcsname\athelp@
{Only certain combinations beginning with @ make sense to me.^^J%
I'll assume you wanted @@ for a printed @.}}{}
\@ifundefined{err@}{\def\err@{\@latexerr}}{}
\@ifundefined{default@tag}%
  {\def\default@tag{%
    \def\tag{\err@{\string\tag\space not allowed here}\@eha}}}
  {}
\@ifundefined{ex@}{\newdimen\ex@}{}
\@ifundefined{minaw@}{\newdimen\minaw@}{}
\@ifundefined{bigaw@}{\newdimen\bigaw@}{}
\newdimen\minCDarrowwidth
\newif\ifCD@
\let\ampersand@\relax
\def\CD{\catcode`\@\active
 \bgroup\relax\let\ampersand@&\iffalse}\fi
 \CD@true\vcenter\bgroup\Let@\restore@math@cr\default@tag
 \tabskip\z@skip\baselineskip20\ex@
 &\hfill$\m@th##$\hfill\crcr}
\def\endCD{\crcr\egroup\egroup\egroup}
\def\CD@check#1#2{\ifCD@\DN@{#2}\else
  \DN@{\err@{@\string#1 not
    allowed outside of the CD environment}\@eha}%
  \fi\next@}
\atdef@>#1>#2>{\ampersand@
  \ifCD@ \global\bigaw@\minCDarrowwidth \else \global\bigaw@\minaw@ \fi
  \setboxz@h{$\m@th\scriptstyle\;{#1}\;\;$}%
  \ifdim\wdz@>\bigaw@\global\bigaw@\wdz@\fi
  \@ifnotempty{#2}{\setbox@ne\hbox{$\m@th\scriptstyle\;{#2}\;\;$}%
    \ifdim\wd@ne>\bigaw@\global\bigaw@\wd@ne\fi}%
 \ifCD@\enskip\fi
   \mathrel{\mathop{\hbox to\bigaw@{\rightarrowfill@\displaystyle}}%
     \limits^{#1}\@ifnotempty{#2}{_{#2}}}%
 \ifCD@\enskip\fi \ampersand@}
\atdef@<#1<#2<{\ampersand@
  \ifCD@ \global\bigaw@\minCDarrowwidth \else \global\bigaw@\minaw@ \fi
  \setboxz@h{$\m@th\scriptstyle\;\;{#1}\;$}%
  \ifdim\wdz@>\bigaw@ \global\bigaw@\wdz@ \fi
  \@ifnotempty{#2}{\setbox@ne\hbox{$\m@th\scriptstyle\;\;{#2}\;$}%
    \ifdim\wd@ne>\bigaw@ \global\bigaw@\wd@ne \fi}%
  \ifCD@\enskip\fi
    \mathrel{\mathop{\hbox to\bigaw@{\leftarrowfill@\displaystyle}}%
      \limits^{#1}\@ifnotempty{#2}{_{#2}}}%
  \ifCD@\enskip\fi \ampersand@}
\begingroup \catcode`\~=\active \lccode`\~=`\@
\lowercase{%
  \global\atdef@)#1)#2){~>#1>#2>}
  \global\atdef@(#1(#2({~<#1<#2<}
}
\endgroup
\atdef@ A#1A#2A{\CD@check{A..A..A}{\llap{$\m@th\vcenter{\hbox
  {$\scriptstyle#1$}}$}\Big\uparrow
  \rlap{$\m@th\vcenter{\hbox{$\scriptstyle#2$}}$}&&}}
\atdef@ V#1V#2V{\CD@check{V..V..V}{\llap{$\m@th\vcenter{\hbox
  {$\scriptstyle#1$}}$}\Big\downarrow
  \rlap{$\m@th\vcenter{\hbox{$\scriptstyle#2$}}$}&&}}
\atdef@={\CD@check={&\enskip\mathrel
  {\vbox{\hrule\@width\minCDarrowwidth\vskip2\ex@\hrule\@width
  \minCDarrowwidth}}\enskip&}}
\atdef@|{\CD@check|{\Big\Vert&&}}
\atdef@\vert{\CD@check\vert{\Big\Vert&&}}
\atdef@.{\CD@check.{&&}}

\makeatletter

\newdimen\paperwidth
\newdimen\paperheight

\def\papersize#1#2{\let\p@persize\relax\paperwidth#1\paperheight#2}

\def\Afour{\papersize{210truemm}{297truemm}}

\let\p@persize\Afour

\let\onesidestyle\@twosidefalse
\let\twosidestyle\@twosidetrue

\def\margins{\@ifnextchar[{\@margins}{\@margins[\z@]}}

\def\@margins[#1]#2#3{
  \p@persize\dimen0 #3\dimen0 .5\dimen0\normalsize%
  \oddsidemargin-1truein\advance\oddsidemargin#2%
  \evensidemargin-1truein\advance\evensidemargin#2%
  \topmargin-1truein\advance\topmargin\dimen0\headsep\dimen0\footskip\dimen0%
  \textwidth\paperwidth\advance\textwidth-#2\advance\textwidth-#2%
  \textheight\paperheight\advance\textheight-#3\advance\textheight-#3%
  \headheight\baselineskip\advance\topmargin-.5\baselineskip%
  \advance\headsep-.5\baselineskip%
  \footheight\baselineskip
  \advance\textwidth-#1\advance\oddsidemargin#1
  \if@twoside\def\@themargin%
    {\ifodd\count\z@\oddsidemargin\else\evensidemargin\fi}\fi}

\def\headlinesep#1{\advance\topmargin\headsep\advance\topmargin -#1
  \advance\topmargin.5\baselineskip\headsep #1\advance\headsep-.5\baselineskip}

\def\footlinesep#1{\normalsize\footskip#1}

\def\headline{\if@twoside\let\n@xt\h@dlin@\else\let\n@xt\h@@dlin@\fi\n@xt}
  
\def\h@dlin@#1#2{%
  \def\@oddhead{%
    {{\leftskip\z@\rightskip\z@\noindent\normalsize#1}}}
  \def\@evenhead{%
    {{\leftskip\z@\rightskip\z@\noindent\normalsize#2}}}}

\def\h@@dlin@#1{%
  \def\@oddhead{{{\leftskip\z@\rightskip\z@\noindent\normalsize#1}}}}

\def\footline{\if@twoside\let\n@xt\f@tlin@\else\let\n@xt\f@@tlin@\fi\n@xt}

\def\f@tlin@#1#2{%
  \def\@oddfoot{%
    {{\leftskip\z@\rightskip\z@\noindent\normalsize#1}}}
  \def\@evenfoot{%
    {{\leftskip\z@\rightskip\z@\noindent\normalsize#2}}}}

\def\f@@tlin@#1{%
  \def\@oddfoot{{{\leftskip\z@\rightskip\z@\noindent\normalsize#1}}}}

\def\normalpage{\global\@specialpagefalse}

\makeatother

\makeatletter

\def\ft{\@ifnextchar[{\ft@s}{\ft@}}
\def\ft@{\ft@@@s[\f@size]}
\def\ft@s[{\@ifnextchar{a}{\ft@sz[}{\ft@@s[}}
\def\ft@@s[{\@ifnextchar{s}{\ft@sz[}{\ft@@@s[}}
\def\ft@@@s[#1]{\ft@sz[at #1pt]}
\def\ft@sz[#1]#2{\font\fonttemp=#2 #1\fonttemp\ignorespaces}

\makeatother

\makeatletter

\input{epsf.sty}

\def\showfig#1#2{\epsfbox{#2}}
\def\fig@#1#2{\leavevmode{\framebox{\figstyl@\strut{ #1 }}}}

\def\figstyle#1{\def\figstyl@{#1}}

\figstyle{\shape{n}\series{m}\selectfont\normalsize}

\def\showfigurestrue{\let\fig\showfig}
\def\showfiguresfalse{\let\fig\fig@}

\makeatother

\showfiguresfalse

\def\smallcircc{\mathop{\mkern3.5mu\hbox{\raise.58ex\hbox{\ft{lcircle10}a}}}}
\def\varemptyset{{\hbox{\raise.21ex\hbox{$\not$}}\mkern.15mu\mathrm{O}\mkern.15mu}}

\let\epsilon\varepsilon
\let\theta\vartheta
\let\phi\varphi
\let\smallcirc\smallcircc
\let\emptyset\varemptyset

\documentstyle[12pt]{article}

\makeatletter

\let\Larg@\Large
\let\hug@\huge

\def\usepackage#1{\input{#1.sty}}

\input{geom.sty}

\let\Large\Larg@
\let\huge\hug@

\def\smallskip{\vskip\smallskipamount}
\def\medskip{\vskip\medskipamount}
\def\bigskip{\vskip\bigskipamount}

\def\mytrivlist{\parsep\parskip\@nmbrlistfalse
  \my@trivlist \labelwidth\z@ \leftmargin\z@
  \itemindent\z@ \def\makelabel##1{##1}}

\def\my@trivlist{\global\@newlisttrue \@outerparskip\parskip}

\def\end#1{\csname end#1\endcsname\@checkend{#1}%
  \expandafter\endgroup\if@endpe\@doendpe\fi
  \if@ignore \global\@ignorefalse \ignorespaces\fi}

\def\partbeforeskip#1{\def\p@rtbeforeskip{#1}}
\def\partstyle#1{\def\p@rtstyl@{#1}}
\def\partdot#1{\def\partd@t{#1}}
\def\partafterskip#1{\def\p@rtafterskip{#1}}
\def\partintrostyle#1{\def\partintr@styl@{#1}}
\def\partintrodot#1{\def\partintr@dot{#1}}
\long\def\partintrosep#1{\long\def\partintr@sep{#1}}
\def\partnewpagetrue{\def\p@rtnewp@ge{\newpage}}
\def\partnewpagefalse{\long\def\p@rtnewp@ge{\par}}

\partbeforeskip{4ex}
\partstyle{\centering\Large\bf}
\partdot{}
\partafterskip{3ex}
\partintrostyle{\large}
\partintrodot{}
\partintrosep{\par}
\partnewpagefalse

\def\partname{Part}
\def\part{\p@rtnewp@ge\addvspace\p@rtbeforeskip\@afterindentfalse\secdef\@part\@spart}

\def\@part[#1]#2{\ifnum \c@secnumdepth >-1\relax  
        \refstepcounter{part}                     
        \def\@tempa{\addcontentsline{toc}{part}}  %
        \expandafter\@tempa\expandafter{\thepart  
          \hspace{1em}#1}\else                    
        \addcontentsline{toc}{part}{#1}\fi        
   {\p@rtstyl@                       
    \ifnum \c@secnumdepth >-1\relax        
      {\partintr@styl@\partname\ \thepart  
       \partintr@dot}\partintr@sep\nobreak 
    \fi                                    
    #2\partd@t\markboth{}{}\par}
    \nobreak                       
    \vskip\p@rtafterskip           
   \@afterheading                  
    }                              

\def\@spart#1{{\p@rtcentering\p@rtstyl@                      
    #1\partd@t\par}                 
    \nobreak                        
    \vskip\p@rtafterskip            
    \@afterheading                  
  }                                 

\def\sectionbeforeskip#1{\def\s@ctbeforeskip{#1}}
\def\sectionstyle#1{\def\s@ctstyl@{#1}}
\def\sectiondot#1{\def\sectiond@t{#1}}
\def\sectionafterskip#1{\def\s@ctafterskip{#1}}
\def\sectionintrostyle#1{\def\sectionintr@styl@{#1}}
\def\sectionintro#1{\def\sectionintr@{#1}}
\def\sectionintrodot#1{\def\sectionintr@dot{#1}}
\def\sectionindenttrue{\def\s@ctind{\parindent}}
\def\sectionindentfalse{\def\s@ctind{\z@}}
\def\sectionafterindenttrue{\def\s@ct@ftind{+}}
\def\sectionafterindentfalse{\def\s@ct@ftind{-}}
\def\sectionafternewlinetrue{\def\s@ct@ftpar{+}}
\def\sectionafternewlinefalse{\def\s@ct@ftpar{-}}

\def\subsectionbeforeskip#1{\def\ss@ctbeforeskip{#1}}
\def\subsectionstyle#1{\def\ss@ctstyl@{#1}}
\def\subsectiondot#1{\def\subsectiond@t{#1}}
\def\subsectionafterskip#1{\def\ss@ctafterskip{#1}}
\def\subsectionintrostyle#1{\def\subsectionintr@styl@{#1}}
\def\subsectionintro#1{\def\subsectionintr@{#1}}
\def\subsectionintrodot#1{\def\subsectionintr@dot{#1}}
\def\subsectionindenttrue{\def\ss@ctind{\parindent}}
\def\subsectionindentfalse{\def\ss@ctind{\z@}}
\def\subsectionafterindenttrue{\def\ss@ct@ftind{+}}
\def\subsectionafterindentfalse{\def\ss@ct@ftind{-}}
\def\subsectionafternewlinetrue{\def\ss@ct@ftpar{+}}
\def\subsectionafternewlinefalse{\def\ss@ct@ftpar{-}}

\def\subsubsectionbeforeskip#1{\def\sss@ctbeforeskip{#1}}
\def\subsubsectionstyle#1{\def\sss@ctstyl@{#1}}
\def\subsubsectiondot#1{\def\subsubsectiond@t{#1}}
\def\subsubsectionafterskip#1{\def\sss@ctafterskip{#1}}
\def\subsubsectionintrostyle#1{\def\subsubsectionintr@styl@{#1}}
\def\subsubsectionintro#1{\def\subsubsectionintr@{#1}}
\def\subsubsectionintrodot#1{\def\subsubsectionintr@dot{#1}}
\def\subsubsectionindenttrue{\def\sss@ctind{\parindent}}
\def\subsubsectionindentfalse{\def\sss@ctind{\z@}}
\def\subsubsectionafterindenttrue{\def\sss@ct@ftind{+}}
\def\subsubsectionafterindentfalse{\def\sss@ct@ftind{-}}
\def\subsubsectionafternewlinetrue{\def\sss@ct@ftpar{+}}
\def\subsubsectionafternewlinefalse{\def\sss@ct@ftpar{-}}

\def\paragraphbeforeskip#1{\def\p@rbeforeskip{#1}}
\def\paragraphstyle#1{\def\p@rstyl@{#1}}
\def\paragraphdot#1{\def\paragraphd@t{#1}}
\def\paragraphafterskip#1{\def\p@rafterskip{#1}}
\def\paragraphintrostyle#1{\def\paragraphintr@styl@{#1}}
\def\paragraphintro#1{\def\paragraphintr@{#1}}
\def\paragraphintrodot#1{\def\paragraphintr@dot{#1}}
\def\paragraphindenttrue{\def\p@rind{\parindent}}
\def\paragraphindentfalse{\def\p@rind{\z@}}
\def\paragraphafterindenttrue{\def\p@r@ftind{+}}
\def\paragraphafterindentfalse{\def\p@r@ftind{-}}
\def\paragraphafternewlinetrue{\def\p@r@ftpar{+}}
\def\paragraphafternewlinefalse{\def\p@r@ftpar{-}}

\def\subparagraphbeforeskip#1{\def\sp@rbeforeskip{#1}}
\def\subparagraphstyle#1{\def\sp@rstyl@{#1}}
\def\subparagraphdot#1{\def\subparagraphd@t{#1}}
\def\subparagraphafterskip#1{\def\sp@rafterskip{#1}}
\def\subparagraphintrostyle#1{\def\subparagraphintr@styl@{#1}}
\def\subparagraphintro#1{\def\subparagraphintr@{#1}}
\def\subparagraphintrodot#1{\def\subparagraphintr@dot{#1}}
\def\subparagraphindenttrue{\def\sp@rind{\parindent}}
\def\subparagraphindentfalse{\def\sp@rind{\z@}}
\def\subparagraphafterindenttrue{\def\sp@r@ftind{+}}
\def\subparagraphafterindentfalse{\def\sp@r@ftind{-}}
\def\subparagraphafternewlinetrue{\def\sp@r@ftpar{+}}
\def\subparagraphafternewlinefalse{\def\sp@r@ftpar{-}}

\sectionbeforeskip{\bigskipamount}
\sectionstyle{\large\bf}
\sectiondot{}
\sectionafterskip{.5\bigskipamount}
\sectionintrostyle{}
\sectionintro{}
\sectionintrodot{.}
\sectionindentfalse
\sectionafterindenttrue
\sectionafternewlinetrue

\subsectionbeforeskip{.8\bigskipamount}
\subsectionstyle{\normalsize\bf}
\subsectiondot{}
\subsectionafterskip{.4\bigskipamount}
\subsectionintrostyle{}
\subsectionintro{}
\subsectionintrodot{.}
\subsectionindentfalse
\subsectionafterindenttrue
\subsectionafternewlinetrue

\subsubsectionbeforeskip{.6\bigskipamount}
\subsubsectionstyle{\normalsize\bf}
\subsubsectiondot{}
\subsubsectionafterskip{.3\bigskipamount}
\subsubsectionintrostyle{}
\subsubsectionintro{}
\subsubsectionintrodot{.}
\subsubsectionindentfalse
\subsubsectionafterindenttrue
\subsubsectionafternewlinetrue

\paragraphbeforeskip{.5\bigskipamount}
\paragraphstyle{\normalsize\bf}
\paragraphdot{.}
\paragraphafterskip{1.25ex}
\paragraphintrostyle{}
\paragraphintro{}
\paragraphintrodot{.}
\paragraphindentfalse
\paragraphafterindenttrue
\paragraphafternewlinefalse

\subparagraphbeforeskip{.5\bigskipamount}
\subparagraphstyle{\normalsize\bf}
\subparagraphdot{.}
\subparagraphafterskip{1.25ex}
\subparagraphintrostyle{}
\subparagraphintro{}
\subparagraphintrodot{.}
\subparagraphindenttrue
\subparagraphafterindenttrue
\subparagraphafternewlinefalse

\def\@startsection#1#2#3#4#5#6{
   \vskip\z@\@tempskipa #4\relax\@afterindenttrue
   \ifdim \@tempskipa <\z@ \@tempskipa -\@tempskipa \@afterindentfalse\fi
   \advance\@tempskipa by\presection
   \if@nobreak \everypar{}\else
     \addpenalty{\@secpenalty}\addvspace{\@tempskipa}%
     \allowbreak\vskip -\presection \fi \@ifstar
     {\@ssect{#1}{#2}{#3}{#4}{#5}{#6}}{\@dblarg{\@sect{#1}{#2}{#3}{#4}{#5}{#6}}}}

\def\@sect#1#2#3#4#5#6[#7]#8{\def\object@type{#1}%
   \ifnum #2>\c@secnumdepth\def\@svsec{}\def\@tempb{}%
      \else\refstepcounter{#1}\def\@svsec{{\csname #1intr@styl@\endcsname%
        {\csname #1intr@\endcsname}}\csname the#1\endcsname%
        \csname #1intr@dot\endcsname\kern1.25ex}%
        \edef\@tempb{\noexpand\numberline{\csname the#1\endcsname}}\fi%
   \@tempskipa #5\relax\def\@tempa{\addcontentsline{toc}{#1}}%
   \ifdim \@tempskipa>\z@%
      \begingroup #6\relax%
        \@hangfrom{\hskip #3\relax\@svsec}{\interlinepenalty \@M #8%
        \csname #1d@t\endcsname\par}%
      \endgroup%
      \csname #1mark\endcsname{#7}%
      \expandafter\@tempa\expandafter{\@tempb #7}%
      \ifautolabel\label*{#8}\fi%
   \else%
      \def\@svsechd{#6\hskip #3\relax%
         \@svsec #8\csname #1mark\endcsname {#7}%
         \expandafter\@tempa\expandafter{\@tempb #7}%
         \ifautolabel\label*{#8}\fi}\fi%
   \@xsect{#5}}

\def\@ssect#1#2#3#4#5#6#7{%
   \ifnum #2>\c@secnumdepth\def\@tempb{}\else \def\@tempb{\numberline{}}\fi%
     \@tempskipa #5\relax\def\@tempa{\addcontentsline{toc}{s#1}}%
     \ifdim \@tempskipa>\z@
        \begingroup #6\relax
           \@hangfrom{\hskip #3}{\interlinepenalty \@M #7%
           \csname #1d@t\endcsname\par}%
        \endgroup
        \csname s#1mark\endcsname{#7}%
        \ifstarredcontents\expandafter\@tempa\expandafter{\@tempb #7}\fi%
        \ifautolabel\label*{#7}\fi%
     \else%
        \def\@svsechd{#6\hskip #3\relax #7\csname s#1mark\endcsname {#7}%
        \ifautolabel\label*{#7}\fi}\fi
   \@xsect{#5}}

\def\section{\@startsection{section}{1}{\s@ctind}
  {\s@ct@ftind\s@ctbeforeskip}{\s@ct@ftpar\s@ctafterskip}{\s@ctstyl@}}
\def\subsection{\@startsection{subsection}{2}{\ss@ctind}
  {\ss@ct@ftind\ss@ctbeforeskip}{\ss@ct@ftpar\ss@ctafterskip}{\ss@ctstyl@}}
\def\subsubsection{\@startsection{subsubsection}{3}{\sss@ctind}
  {\sss@ct@ftind\sss@ctbeforeskip}{\sss@ct@ftpar\sss@ctafterskip}{\sss@ctstyl@}}
\def\paragraph{\@startsection{paragraph}{4}{\p@rind}
  {\p@r@ftind\p@rbeforeskip}{\p@r@ftpar\p@rafterskip}{\p@rstyl@}}
\def\subparagraph{\@startsection{subparagraph}{4}{\sp@rind}
  {\sp@r@ftind\sp@rbeforeskip}{\sp@r@ftpar\sp@rafterskip}{\sp@rstyl@}}

\def\statementabove#1{\def\th@bove{#1}}
\def\statementstyle#1{\def\thstyl@{#1}}
\def\statementbelow#1{\def\thb@low{#1}}
\def\statementindentfalse{\let\thind@nt\relax}
\def\statementindenttrue{\let\thind@nt\indent}

\def\statementintrostyle#1{\def\thintr@style{#1}}
\def\statementintrodot#1{\def\thintr@dot{#1}}
\def\statementintrosep#1{\def\thintr@sep{#1}}
\def\statementintrobrackets#1#2{\def\thintr@left{#1}\def\thintr@right{#2}}

\statementabove{\medskip}
\statementstyle{\sl}
\statementbelow{\medskip}
\statementindenttrue

\statementintrostyle{\normalshape\bf}
\statementintrodot{.}
\statementintrosep{\kern1.25ex}
\statementintrobrackets{(}{)}

\def\@thskip{\dimen0\lastskip\vskip-\dimen0%
  \th@bove\dimen1\lastskip\vskip-\dimen1%
  \ifdim\dimen0>\dimen1\else\dimen0\dimen1\fi\vskip\dimen0}

\long\def\@@newtheorem#1#2#3{%
  \newenvironment{#3}%
    {\def\object@type{#3}\par\@thskip#1%
     \@ifnextchar[{\@enva{#3}{\thstyl@{#2}}}{\@envb{#3}{\thstyl@{#2}}}}%
    {\end{#3@}}%
  \@ifnextchar[{\@othm{#3@}}{\@nnthm{#3}}}

\def\theoremintro#1{\thintr@style{#1}\thintr@dot\thintr@sep}

\def\thrmintro#1#2{#1%
     \def\@tempa{#1}\ifx\@tempa\@empty\else
     \expandafter\let\expandafter\@tempa#2
     \ifx\@tempa\relax\else\kern1ex\fi\fi#2}

\def\@begintheorem#1#2{%
        \mytrivlist\item[\thind@nt\hskip\labelsep
        {\theoremintro{\thrmintro{#1}{#2}}}\hskip-\labelsep]}

\def\@opargbegintheorem#1#2#3{%
        \mytrivlist\item[\thind@nt\hskip\labelsep
        \theoremintro{\thrmintro{#1}{#2}\kern1ex
        \thintr@left{#3}\thintr@right}\hskip-\labelsep]%
        \ifautolabel\label*{#3}\fi}

\def\@endtheorem{\endtrivlist\thb@low}

\def\proofabove#1{\def\pf@bove{#1}}
\def\proofstyle#1{\def\pfstyl@{#1}}
\def\proofbelow#1{\def\pfb@low{#1}}
\def\proofindentfalse{\let\pfind@nt\relax}
\def\proofindenttrue{\let\pfind@nt\indent}

\def\proofintrostyle#1{\def\pfintr@style{#1}}
\def\proofintrodot#1{\def\pfintr@dot{#1}}
\def\proofintrosep#1{\def\pfintr@sep{#1}}
\def\proofintrobrackets#1#2{\def\pfintr@left{#1}\def\pfintr@right{#2}}

\proofabove{\medskip}
\proofstyle{}
\proofbelow{\medskip}
\proofindenttrue

\proofintrostyle{\sl}
\proofintrodot{.}
\proofintrosep{\kern1.25ex}
\proofintrobrackets{of\kern1ex}{}

\def\@pfskip{\dimen0\lastskip\vskip-\dimen0%
  \pf@bove\dimen1\lastskip\vskip-\dimen1%
  \ifdim\dimen0>\dimen1\else\dimen0\dimen1\fi\vskip\dimen0}

\renewenvironment{proof}%
  {\@pfskip\mytrivlist\item[\pfind@nt]\@ifnextchar[{\pro@f}{\pro@f[\prooftag]}}
  {\ifvoid\provedbox\else\hproved\fi\endtrivlist\pfb@low}

\def\pro@f[#1]{\setbox\provedbox\hbox{\provedboxcontents{#1}}\proofintro{#1}}

\def\proofintro#1{\expandafter\def\expandafter\@tempa\expandafter{#1}%
  {\pfintr@style{Proof\ifx\@tempa\empty\else\kern1ex\pfintr@left{#1}%
  \pfintr@right\fi}\pfintr@dot\pfintr@sep}\pfstyl@\ignorespaces}

\def\provedmark#1{\def\prm@rk{#1}}
\def\provedsep#1{\def\prs@p{#1}}

\provedmark{$\square$}
\provedsep{\kern1.25ex}

\def\provedtexttrue{\def\prb@x##1{\fbox{\small##1}}}
\def\provedtextfalse{\def\prb@x##1{\prm@rk}}
\def\provedmarkrighttrue{\let\prhf@l\hfill}
\def\provedmarkrightfalse{\let\prhf@l\relax}

\provedtextfalse
\provedmarkrightfalse

\def\provedboxcontents#1{\expandafter\def\expandafter\@tempa\expandafter{#1}%
  \ifx\@tempa\empty\prm@rk\else\prb@x{#1}\fi}

\def\proved{\ifmmode\eqno{\box\provedbox}\else\hproved\fi}

\def\hproved{\unskip\nobreak\prhf@l\penalty50\prs@p\hbox{}\nobreak\prhf@l
  \box\provedbox{\finalhyphendemerits=0\par}}

\def\captionstyle#1{\def\c@ptstyl@{#1}}
\def\captionintrostyle#1{\def\c@pintr@style{#1}}
\def\captionintrodot#1{\def\c@pintr@dot{#1}}
\def\captionintrosep#1{\def\c@pintr@sep{#1}}

\captionstyle{\small\sf}
\captionintrostyle{\bf}
\captionintrodot{.}
\captionintrosep{\hskip1.25ex}

\long\def\@makecaption#1#2{%
    \vskip\captionskip
    \setbox\@tempboxa\hbox{%
      \ifproofing\@ifundefined{the@label}{}
        {\hbox to 0pt{\vbox to 0pt{\vss\hbox{\tiny\the@label}\bigskip}\hss}}\fi
      \c@ptstyl@{\c@pintr@style #1\c@pintr@dot}\ignorespaces #2}%
    \@captionwidth=\hsize \advance\@captionwidth-2\@captionmargin
    \ifdim \wd\@tempboxa >\@captionwidth {%
        \rightskip=\@captionmargin\leftskip=\@captionmargin
        \unhbox\@tempboxa\par}%
      \else
        \hbox to\hsize{\hfil\box\@tempboxa\hfil}%
    \fi}

\def\end@Float#1{%
  \expandafter\caption\expandafter[\the@title]{%
   {\c@pintr@style%
   \ifx\the@caption\empty\ifx\the@title\empty\else\c@pintr@sep\fi\else\c@pintr@sep\fi
    \the@title\ifx\the@caption\empty\else\ifx\the@title\empty\else
    \c@pintr@dot\c@pintr@sep\fi\fi}%
   \ignorespaces\the@caption%
  \expandafter\label\expandafter*\expandafter{\the@label}}%
  \end{#1}}

\@definecounter{bibenumi}

\def\thebibliography#1{\section*{\refname}%
 \list{[\arabic{bibenumi}]}{\settowidth\labelwidth{[#1]}%
 \leftmargin\labelwidth\advance\leftmargin\labelsep\usecounter{bibenumi}}%
 \def\newblock{\hskip .11em plus .33em minus .07em}%
 \sloppy\clubpenalty4000\widowpenalty4000\sfcode`\.=1000\relax}

\proofingfalse
\autolabelfalse

\showfiguresfalse

\newtheorem{stat}{\statname}  \unnumbered{stat}

\newtheorem{nstat}{\nstatname}[section]

\newtheorem{lemma}[nstat]{Lemma}
\newtheorem{proposition}[nstat]{Proposition}
\newtheorem{theorem}[nstat]{Theorem}
\newtheorem{corollary}[nstat]{Corollary}

\newtheorem{remark}[nstat]{Remark}

\papersize{215truemm}{275truemm}
\margins{3cm}{2.8cm}
\footlinesep{1.4cm}
\headline{\hfill}
\footline{\small\hfill--\kern1ex\thepage\kern1ex--\hfill}
\makeatletter\c@totalnumber8\c@topnumber8\makeatother
\showfigurestrue

\sectionbeforeskip{1.5\bigskipamount}
\sectionstyle{\centering\normalsize\bf}
\sectionafterskip{\bigskipamount}
\sectionafterindenttrue

\paragraphbeforeskip{\bigskipamount}
\paragraphstyle{\centering\normalsize\sl}
\paragraphafterskip{.5\bigskipamount}
\paragraphafternewlinetrue
\paragraphafterindenttrue
\paragraphdot{}

\statementintrostyle{\sc}

\captionstyle{\small}
\captionintrostyle{\sc}

\newcommand{\R}{R} 
\newcommand{\C}{C} 
\newcommand{\B}{B} 
\renewcommand{\S}{S}

\newcommand{\Cl}{\mathop{\mathrm{Cl}}\nolimits} 
\newcommand{\Int}{\mathop{\mathrm{Int}}\nolimits} 
\newcommand{\Bd}{\mathop{\mathrm{Bd}}\nolimits}
\newcommand{\Map}{\mathop{\mathrm{Map}}\nolimits}
\newcommand{\id}{\mathop{\mathrm{id}}\nolimits}

\begin{document}

\title{\large\bf COMPACT STEIN SURFACES WITH BOUNDARY\\
AS BRANCHED COVERS OF $\B^4$}
\author{\sc\normalsize A. Loi\\
\sl\normalsize Dipartimento di Matematica e Fisica\\[-3pt]
\sl\normalsize Universit\`a di Sassari -- Italia\\
\tt\small loi@ssmain.uniss.it
\and 
\sc\normalsize R. Piergallini\\
\sl\normalsize Dipartimento di Matematica e Fisica\\[-3pt]
\sl\normalsize Universit\`a di Camerino -- Italia\\
\tt\small pierg@camserv.unicam.it}
\date{}

\maketitle

\begin{abstract}
\baselineskip13.5pt
\smallskip
\noindent
We prove that Stein surfaces with boundary coincide up to orientation preserving
diffeomorphisms with simple branched coverings of $\B^4$ whose branch set is a positive
braided surface. As a consequence, we have that a smooth oriented $3$-manifold is
Stein fillable iff it has a positive open-book decomposition.



\end{abstract}

\section*{Introduction}

Compact Stein surfaces with (stricly pseudoconvex) boundary play an important role for
the contact topology of 3-manifolds, due to the fact that their boundaries carry
natural tight contact structures, given by the complex tangenties.

It is worth remarking that, this is one of the only two known general ways for
producing tight contact structures, the other one being perturbation of taut
foliations (cf. \cite{ET98}). On the other hand, Stein surfaces with boundary can also
be used to define invariants for fillable contact structures (see \cite{G98} and
\cite{LM98}).

A topological characterization of compact Stein surfaces has been given by Eliashberg
in term of handle decomposition, by using the notion of Legendrian surgery (cf.
\cite{E90} and \cite{G98}). In \cite{G98}, Gompf developed a Legendrian version of the
Kirby calculus on framed links, in order to construct and study fillable contact
3-manifolds. In the same paper, he conjectured that the Poincar\'e homology sphere with
reversed orientation could not be Stein fillable.
This conjecture has been proved in \cite{L98} by Lisca. Successively, Ethnyre and
Honda showed that the Poincar\'e homology sphere with reversed orientation cannot
carry any tight contact structure (see \cite{EH99}). However, we still have no
general way for establishing whether a given 3-manifold has such a contact
structure or not.

In this paper we propose an alternative approach to the topology of Stein surfaces with
boundary, representing them as branched covers of $\B^4$.
Namely, starting with a Legendrian handle decomposition of $X$, the lifting surgery
method introduced by Montesinos in \cite{M78}, gives us a covering $p:X \to \B^4$,
whose branch set is a non-singular ribbon (real) surfaces $S \subset \B^4$. Then, we
can apply the Rudolph's braiding process to $S$ (cf. \cite{R83a}) in order to make $S$
into a braided surface in $\B^2 \times \B^2 \cong \B^4$. The crucial point is that,
performing all the operations in the proper way, the resulting braided surface is
positive. By \cite{R83}, this means that we can assume $S$ to be analytic. At this
point, the Grauert-Remmert theory of analitically branched coverings (see \cite{DG94}
or \cite{GR58}), allows us to conclude that $p$ itself can be assumed analytic up to
orientation preserving diffeomorphisms. Viceversa, it is not difficult to prove that
any analytical branched cover of $\B^4$ is orientation preserving diffeomorphic to a
Stein surface with boundary.

By composing the branched covering $p$ with the projection $\B^4 \cong \B^2 \times
\B^2 \to \B^2$, we get a positive Lefschetz fibration $f: X \to \B^2$. In fact, under
some natural restrictions, any Lefschetz fibration over $\B^2$ factors in such a way.
This gives us a further topological characterization of the compact Stein surfaces with
boundary as positive Lefschetz fibrations of $\B^2$.
Looking at the boundary, we immediately get a corresponding fillability
criterion in terms of positive open books.

\medskip
The paper is organized as follows. In section \ref{fib/sec} we prove some preliminary
results relating Leschetz fibrations with covering branched over braided surfaces.
Section \ref{stein/sec} is interely devoted to prove our main theorem, that is the
characterizations of compact Stein surfaces with boundary as branched covering of
$\B^4$ and as Lefschetz fibrations over $\B^2$ (theorem \ref{stein/thm}). Finally, in 
section \ref{fil/sec} we use this characterization in order to obtain the above mentioned
fillability criterion (theorem \ref{fil/thm}). 


\section{Lefschetz fibrations\label{fib/sec}}

Let $X$ be a smooth oriented connected compact $4$-manifold with (possibly empty)
boundary and $Y$ be a smooth oriented connected compact surface with (possibly empty) 
boundary. A smooth map $f:X \to Y$ is called a {\sl Lefschetz fibration} over $Y$ iff
the following properties hold:
\begin{enumerate}
\item \vskip-\topsep\smallskip
$f$ has finitely many singular values $y_1,\dots,y_n \in \Int Y$ (the {\sl branch
points} of $f$) and the restriction of $f$ over $Y - \{y_1,\dots,y_n\}$ is a
locally trivial fiber bundle whose fiber $F$ is an oriented compact surface
with (possibly empty) boundary (the {\sl regular fiber} of $f$);
\item for each $i=1,\dots,n$, there is only one singular point $x_i \in \Int X$ over
the branch point $y_i$ and the monodromy of a counterclockwise meridian loop around
$y_i$ is given by $\delta_i^{\epsilon_i}$, where $\delta_i$ is the right-handed Dehn
twist along $d_i \subset \Int F$ and $\epsilon_i = \pm 1$ ($x_i$ is
called {\sl positive} or {\sl negative} depending on $\epsilon_i$).
\end{enumerate}\vglue-\medskipamount

We say that $f$ is {\sl positive} iff all its singular points $x_i$ are positive and
that $f$ is {\sl allowable} iff all the loops $d_i$ are homologically non-trivial in
$F$.

\medskip
A Lefschetz fibration $f:X \to Y$ is completely determined, up to
orientation preserving diffeomorphisms, by the branch points $y_1, \dots, y_n \in \Int
Y$ and by its restriction over $Y -\{y_1,\dots,y_n\}$. On the other hand, any locally
trivial fiber bundle over $Y - \{y_1,\dots,y_n\}$ satisfying (a) and (b)
uniquely extends to a Lefschetz fibration.
In fact, the structure of $f$ over a small disk $D_i$ centered at $y_i$ is given by
the following commutative diagram, where: $T(\delta_i^{\epsilon_i})$ is the mapping
torus of $\delta_i^{\epsilon_i}$ and $\pi:T(\delta_i^{\epsilon_i}) \to \S^1$ is the
canonical projection; the singular fiber $F_{y_i} \cong F / d_i$ has a transversal
self-intersection at $x_i$, which is positive or negative depending on
$\epsilon_i$;\break
$h$ and $k$ are orientation preserving diffeomorphisms such that, denoting with\break
$i_{s,t}: F \to T(\delta_i^{\epsilon_i}) \times (0,1]$ the canonical inclusion defined
by
$i_{s,t}(x) = ([x,s],t)$ and putting $k_{s,t} = k\smallcirc i_{s,t}: F \to F_{h(s,t)}
\subset f^{-1}(D_i - \{y_i\})$, we have $k_{s,t}(d_i) \to x_i$ as $t \to 0$.

$$\begin{CD}
T(\delta_i^{\epsilon_i}) \times (0,1]@>k>> f^{-1}(D_i - \{y_i\}) 
&\ \ \subset \ \ &f^{-1}(D_i) &\ \ \supset \ \ & F_{y_i}\cr
@VV{\pi\times\hbox{id}}V @VVV @VVf_{|}V @VVV\cr
\S^1 \times (0,1] @>h>> D_i - \{y_i\} &\ \ \subset \ \ &D_i & \supset & \{y_i\}
\end{CD}$$

\medskip
For any $i=1,\dots,n$, there are local complex coordinates $(z_1,z_2)$
of $X$ and $z$ of $Y$, respectively centered at $x_i$ and at $y_i$, making $f$ into
the complex map $(z_1,z_2) \mapsto z = z_1^2 + z_2^2$. Moreover, such coordinates can
be chosen orientation preserving iff $x_i$ is a positive singular point.
In other words, $f$ is locally a complex Morse function. This fact could be used to get
a natural handle decomposition of $Y$. For a detailed discussion of the topology of
Leschetz fibrations we refer to \cite{GS99}.

\medskip
If $\Bd Y\neq \emptyset$, the observations above say that a Lefschetz fibration $f:X
\to Y$ is uniquely determined, up to orientation preserving diffeomorphisms, by its
monodromy
$\phi_f:\pi_1(Y -
\{y_1,\dots,y_n\},\ast) \to \Map F$ and that $\phi_f$ can be an arbitrary
homomorphism satisfying the property (b).

For $Y = \B^2$, the monodromy $\phi_f$ can be represented by an arbitrary
sequence of Dehn twists $\delta_1^{\epsilon_1}, \dots, \delta_n^{\epsilon_n}$ along
simple loops $d_1, \dots, d_n \subset \Int \B^2$, giving the monodromies of
counterclockwise meridian loops around the branch points $y_1,\dots,y_n$, which freely
generate
$\pi_1(\B^2 - \{y_1,\dots,y_n\},\ast)$.

\medskip
In order to describe Lefschetz fibrations in terms of branched coverings, we introduce
the notion of braided surface in a product of surfaces (cf. \cite{R83a} for the case of
$\B^2 \times \B^2$).

\medskip
Let $Y$ and $Z$ be smooth oriented connected compact
surfaces. A regularly embedded smooth compact surface $S \subset Y \times Z$
is a {\sl braided surface} over $Y$ iff the restriction of the canonical projection
$\pi_{Y|S}:S \to Y$ is a simple branched covering. 

We observe that $S$ is oriented as branched cover of $Y$ and $\Bd S$ is an oriented
link in $\Bd Y \times Z$ which intersects $C \times Z$ in a closed braid, for every
component $C$ of $\Bd Y$. Furthermore, $\pi_{Y|S}$ has finitely many singular value
$y_1, \dots, y_n \in \Int Y$ and over each $y_i$ there is only one singular point $s_i
\in \Int S$ for $\pi_{Y|S}$. We call $s_1, \dots, s_n$ the {\sl twist points} of $S$. 

For any twist point $s_i$ of $S$, there are fiber preserving local complex coordinates
$(w,z)$ of $Y \times Z$ centered at $s_i$ making $S$ into the surface $w = z^2$.
We say that $s_i$ is a positive twist point iff such coordinates can be choosen
orientation preserving (with respect to the product orientation of $Y \times Z$)
and a negative twist point otherwise. We call $S$ a {\sl positive} braided surface iff
all its twist points are positive.

\medskip
The following theorem on positive braided surfaces in $\B^2 \times \B^2$ will
be used in the next section. Its proof is implicit in \cite{R83} (see remark 4.4 in
\cite{R83a} and observe that any positive braided surfaces in $\B^2 \times \B^2$ has a
quasipositive band presentation).

\begin{theorem}[Rudolph] \label{rud/thm}
A braided surface $S \subset \B^2 \times \B^2$ is positive iff it is isotopic to
the intersection of a complex analytic curve with $\B^2\times \B^2 \subset \C^2$.
\end{theorem}

Now, we come to the relation between Lefschetz fibrations with fiber $F$ over a surface
$Y$ and branched coverings of products $Y \times Z$ (typically $Z \cong \S^2$ for $F$
closed and $Z \cong \B^2$ for $F$ bounded) with branch surfaces $S \subset Y
\times Z$ braided over $Y$.

\begin{proposition} \label{bclf/thm}
Let $Y$ and $Z$ be smooth oriented connected compact surfaces and let $p:X \to Y\times
Z$ be a simple branched covering whose branch set is a surface $S \subset Y
\times Z$ braided over $Y$. Then, the composition $f = {\pi_Y \smallcirc p}: X \to Y$
is a Lefschetz fibration which has the same branch points of $\pi_{Y|S}$ and one
positive (resp. negative) singular point over each positive (resp. negative) twist
point of $S$. Moreover, if $\,\Bd Z \neq \emptyset$ then the regular fiber of $f$ has
no closed components and $f$ is allowable.
\end{proposition}

\begin{proof}
Of course, $f$ is regular at each regular point of $p$. Furthermore, given $x
\in X$ singular point of $p$, we have $p(x) \in S$ and 
$T_xf(T_xX) = T_{p(x)}\pi_Y(T_xp(T_xX)) = T_{p(x)}\pi_Y(T_{p(x)}S)$, hence $x$ a
singular point of $f$ iff $p(x)$ is a twist point of $S$.

Now, let $s_1,\dots,s_n \in S$ the twist points of $S$ and $y_1, \dots, y_n \in Y$
their projections by $\pi_Y$. Then, $f$ is regular over $Y - \{y_1,\dots,y_n\}$ and, by
compactness, it satisfies property (a) of Lefschetz fibrations, the regular fiber
$F \cong f^{-1}(y)$ with $y \neq y_1,\dots,y_n$ being simple covering of $Z \cong \{y\}
\times Z$ branched over the (transversal) intersection with $S$, by the restriction of
$p$.

On the other hand, since $p$ is simple, over each singular value $y_i$ there is only
one singular point $x_i$. In order to verify property (b) of Lefschetz fibrations,
we have to check that the monodromy around each $y_i$ is a Dehn twist.

Let $(w,z)$ be local fiber preserving complex coordinates of $Y \times Z$ centered at
$s_i$ and making $S$ into the surface $w=z^2$. We can assume that $w$ is orientation
preserving on $Y$, so that $t \mapsto w(t) = \rho e^{2\pi i t}$, with $\rho > 0$
sufficiently small, is a counterclokwise parametrization of a simple loop $l_i
\subset Y$ around $y_i$. 

Then $S \cap (l_i \times Z)$ is the closed braid in $l_i \times Z$,
corresponding to a half twist around an arc $a \subset \{w(0)\} \times \Int Z$ between
two branch points of the restriction of $p$ over $\{w(0)\} \times Z$, whose meridians
have the same monodromy. Such a half twist is right-handed (resp. left-handed) if
$s_i$ is a positive (resp. negative) twist point of $S$ and lifts to the right-handed
(resp. left-handed) Dehn twist along the unique simple loop $d$ contained in
$p^{-1}(a) \subset \Int f^{-1}(w(0)) \cong \Int F$ (cf. \cite{BE79}, lemma 4.2), which
represents the monodromy of $l_i$.

Finally, assuming $\Bd Z \neq \emptyset$, we have that each component of the regular
fiber $F$ has non-empty boundary, since it is a branched covering of $Z$. Similarly, 
for the loop $d \subset F$ considered above, we have that each component of $F - d$ has
non-empty boundary. Then, we can conclude that $f$ is allowable if $\Bd Z \neq
\emptyset$.
\end{proof}

The following proposition shows that any allowable Lefschetz fibration over $Y$ whose
fiber is connected with (possibly empty) connected boundary, can be obtained as in
proposition \ref{bclf/thm} from a quite special branched covering if
$\Bd Y \neq \emptyset$.

\begin{proposition} \label{lfbc/thm}
Let $f:X \to Y$ be an allowable Lefschetz fibration with regular fiber $F$. If $F$ and
$\,\Bd F$ are connected and $\,\Bd Y \neq \emptyset$, there exists a 3-fold simple
branched covering $p:X \to Y \times Z$ whose branch set is a surface $S \subset Y
\times Z$ braided over $Y$, with $Z \cong \S^2$ if $F$ is closed and $Z
\cong \B^2$ otherwise, such that $f = \pi_Y \smallcirc p$.
\end{proposition}

\begin{proof}First of all, since $F$ and $\Bd F$ are connected, there exists a 3-fold
simple branched covering $q:F \to Z$, with $Z$ as in the statement, such that any
Dehn twist of $F$ along a non-separating simple loop can be realized, up to isotopy, as
the lifting of a half twist around an arc in $Z$ between two branch points of $q$,
whose meridians have the same monodromy (see \cite{BW85} and remember that all the
non-separating simple loops in $F$ are equivalent).
Then, any element of $\Map F$ can be represented by the lifting of a diffeomorphism of
$Z$ onto itself isotopic to the identity, since Dehn twists along non-separating simple
loops generate $\Map F$.

Let $y_1, \dots, y_n \in \Int Y$ be the branch points of $f$ and $A_1, \dots, A_n
\subset Y$ be disjoint disks such that $y_i \in \Int A_i$ and $A_i \cap \Bd Y$ is an
arc in $\Bd A_i$ for every $i = 1, \dots, n$. 
Then, the restriction of $f$ over $Y_0 = \Cl(Y - (A_1 \cup \dots \cup A_n))$ is a
locally trivial fiber bundle. 

Given a band presentation $Y_0 \cong \B^2 \cup H_1 \dots \cup H_m$ with bands
($=$ 1-handles)
$H_1,
\dots, H_m$, we construct a branched covering $p_0:X_0 \to Y_0 \times Z$ as follows:
start with the covering
$\id_Y
\times\,q: Y_0
\times F \to Y_0
\times Z$; cut each $H_j \times F$ along $t_j \times F$ and each
$H_j \times Z$ along $t_j \times Z$, where $t_j$ is a transversal arc for the band
$H_j$; glue them back respectively by $\id_{t_j}\! \times\, \phi_f(e_j)$ and
$\id_{t_j}\!\times\, h_j$, where $\phi_f(e_j) \in \Map F$ is the monodromy of a
simple loop $e_j$ which goes once through $H_j$ and $h_j: Z \to Z$ is a
homeomorphism isotopic to the identity which lifts to $\phi_f(e_i)$ by means of $q$.

In order to extend $p_0$ to a branched covering $p: X \to Y$, we consider a branched
covering $r:W \to \B^2 \times Z$ whose branch set is a surface braided over $\B^2$ with
only one positive twist point over $0$ and whose restriction over $\S^1_- \times Z$
coincides with $\id_{\S^1_-}\!\times q$. As we have seen in the proof of proposition
\ref{bclf/thm}, the composition $\pi_{\B^2} \smallcirc r$ is a Lefschetz fibration
branched over $0$ with regular fiber $F$, such that the monodromy of a 
counterclockwise meridian loop around $0$ is a right-handed Dehn twist along a
non-separating simple loop $\delta \subset \Int F$.

Now, for any $i = 1, \dots, n$, let $\phi_f(l_i) = \delta_i^{\epsilon_i}$, where $l_i
\subset A_i$ is a counterclockwise meridian loop around $y_i$, $\delta_i$ is the
right-handed Dehn twist along $d_i \subset \Int F$ and\break  $\epsilon_i = \pm 1$.
Since $f$ is allowable, $d_i$ cannot separate $F$, so there exist diffeomorphisms $k_i
= k_i' \times k_i'':\S^1_- \times Z \to A_i \times Z$ and  $\widetilde k_i = k_i'
\times \widetilde k_i'':\S^1_- \times F \to A_i \times F$ such that:\break
$k_i'$ preserve or invert the orientation according to $\epsilon_i$;
$k_i''$ is orientation preserving and lifts to $\widetilde k_i''$ with respect to $q$; 
$\widetilde k_i''(d) = d_i$.
Then, assuming that the arcs $a_i = A_i \cap Y_0$ do not meet the 1-handles $H_j$, we
can glue $n$ copies of $r$ to $p_0$, by means of the diffeomorphisms $k_i$ and
$\widetilde k_i$.

Calling $p$ the branched covering of $Y$ obtained in this way, we have that $\pi_Y
\smallcirc p$ is a Lefschetz fibration whose branch points a monodromy coincide with
that ones of $f$, by proposition \ref{bclf/thm} and its proof. So, up to orientation 
preserving diffeomorphisms, $\pi_Y \smallcirc p = f$ and in particular the total 
space of $p$ is $X$.
\end{proof}

\begin{remark} \label{lfbc/rem} 
Proposition \ref{lfbc/thm} does not hold in general if $\,\Bd Y = \emptyset$ (see \cite{F99}
for hyperelliptic Lefschetz fibrations). In fact, to deal with this case, we should
allow the braided surface $S$ to have node and cusp singularities (cf. \cite{P95}). 
The connection requirement for $F$ and $\Bd F$ could perhaps be removed, by considering
branched coverings of order greater that 3.
\end{remark}

We conclude this section by observing that, for a Lefschetz fibration $f:X \to \B^2$,
the condition of having connected fiber with connected boundary, does not imply any
restriction on the total space $X$. This fact will be needed in the next section.

\begin{proposition} \label{lfB2/thm}
If $f:X \to \B^2$ is a Lefschetz fibration over $\B^2$, then the regular fiber of $f$
is connected and there exists a Lefschetz fibration $g:X \to \B^2$ whose fiber has
connected boundary. Moreover, for $f$ allowable and/or positive, we can take
$g$ allowable and/or positive as well.
\end{proposition}

\begin{proof}
The connection of $F$ follows immediately from the connection of $X$, since the
monodromy of $f$ is generated by Dehn twists, so it preserves the components of $F$.
We also observe that, for the same reason, the monodromy of $f$ fixes the boundary of
$F$.

Now, if $\Bd F = \emptyset$ or $\Bd F$ is already connected, we can set $g = f$.
Othewise, in order to connect the boundary of $F$, we consider the following plumbing
operation for Lefschetz fibrations with connected bounded fiber, which is analogous to
the operation (A) introduced by Harer in \cite{H82} for open-book decomposition.

Let $F' = F \cup H$ the surface obtained by gluing an oriented band $H$ to $F$ (we are
assuming $\Bd F \neq \emptyset$) and $d \subset \Int F'$ be a simple loop which goes
once through $H$ (we are also assuming $F$ connected).
Then, we consider the new Lefschetz fibration $f':X' \to Y$ with regular fiber
$F'$, branch points $y_1, \dots, y_n,y_{n+1}\in\Int \B^2$ and respective monodromies
$\delta_1^{\epsilon_1}, \dots, \delta_n^{\epsilon_n},\delta$, where $y_1, \dots, y_n$
are the branch points of $f$, $\delta_1^{\epsilon_1}, \dots, \delta_n^{\epsilon_n}$ are
the respective monodromies for $f$ thought as Dehn twists of $F'$ and $\delta$ is
the right-handed Dehn twist along $d$. 

By the definition of $f'$, we get $X' \cong
X$, in fact $X'$ can be obtained by adding to $X$ a cancelling pair of handles: one
1-handle $\B^2 \times H$ glued to $\B^2 \times \Bd F \subset \Bd X$ (remember that
the monodromy of $f$ fixes $\Bd F$), due to the change of the fiber, and one 2-handle
attached along $\{s\} \times d \subset \{s\} \times G \subset \Bd(X \cup (\B^2 \times
H))$ with $s \in \S^1$, due to the new branch point $y_{n+1}$ (cf. \cite{F99} and
\cite{K80}).
On the other hand, if $\Bd F$ is not connected and the band $H$ joins two different
components of $\Bd F$, then $\Bd F'$ has one component less than $\Bd F$ and $d$ is
non-separating in $F'$. 

Then we can get the required Lefschetz fibration $g$ from $f$, by iterating the
plumbing operation, until the boundary of the fiber becomes connected.
\end{proof}

\begin{remark} \label{plumbing/rem}
For a Lefschetz fibration $f = \pi_{\B^2} \smallcirc p$, with $p:X \to \B^2 \times
\B^2$ simple covering branched over a braided surface $S \subset \B^2 \times \B^2$, a
plumbing operation on $f$ corresponds to a stabilization of $S$, consisting in the
addition of one sheet connected to $S$ by means of one positive twist point.
\end{remark}


\section{Stein surfaces\label{stein/sec}}

We recall that, a smooth real-valued function $f:X \to R$ on a complex manifold $X$ is
called {\sl plurisubharmonic} (resp. {\sl strictly plurisubharmonic}) iff the complex
Hessian $H f = (\partial^2 f/\partial z_i \partial \bar z_j)$ is everywhere
positive semidefinite (resp. definite) for any local complex coordinates $(z_1, \dots,
z_n)$. Of course, both these properties are invariant under biholomorphisms of $X$.
Moreover, plurisubharmonicity (but not strict plurisubharmonicity) is preserved under
composition with holomorphic functions on the right and with non-decreasing convex
funtions on the left (see \cite{GR65} or \cite{P94}).

\medskip
A {\sl Stein surface} is a non-singular complex surface $X$ which 
admits a proper strictly plurisubharmonic function $f:X \to [0,+\infty)$ such
that $\Bd X$ is a level set.

If $X \subset \C^n$ is a non-singular complex surface properly embedded in $\C^n$, then
the restriction to $X$ of the function $z \mapsto |z|^2$ is a proper strictly
plurisubharmonic function, hence $X$ is a Stein surface.  In this way we get all the
Stein surfaces without boundary, up to biholomorphisms, since any Stein surface
without boundary can be properly holomorphically embedded in some $\C^n$ (see
\cite{GR58} or \cite{GR65}).

If $X$ is a Stein surface without boundary and $f:X \to [0,+\infty)$ is a proper
strictly plurisubarmonic function, then the sublevel set $f^{-1}([0,c])$ is a compact
Stein surface with boundary $f^{-1}(c)$, for any regular value $c > 0$. Any
compact Stein surface has non-empty boundary and can be embedded in a Stein surface
without boundary as a sublevel set of some proper plurisubharmonic function as above.

\medskip
Any Stein surface $X$ has a (possibly infinite) handle decomposition, induced by a
plurisubharmonic Morse function, with handles of indices $\leq 2$ (see \cite{M63}).

In particular, for $X$ compact we get $X \cong X_1 \cup H_1 \cup \dots \cup H_m$,
where $X_1$ is obtained by attaching 1-handles to $\B^4$ and the $H_i$'s are
2-handle attached to $X_1$.\break
It turns that the $H_i$'s are attached to $X_1$ in a quite
special way. In fact, the attaching knot $K_i \subset \Bd X_1$ of each 2-handle
$H_i$ is Legendrian with respect to the standard contact structure of $\Bd X_1
\cong \#_{n\,}\S^1 \times \S^2$ and the attaching framing is the Legendrian framing 
of $K_i$ with one left-handed twist added (see \cite{G98} or \cite{GS99} for more
details).

We call {\sl Legendrian} such a 2-handle $H_i$. For our aims, it will suffice to know
how to represent Legendrian 2-handles in terms of framed links. The translation in the
language of framed links is widely discussed in \cite{G98} and \cite{GS99}, so we
limit ourselves to describe the final form of the resulting framed link.

\begin{Figure}[htb]{}{}{}
\centerline{\fig{Front projection without 1-handles}{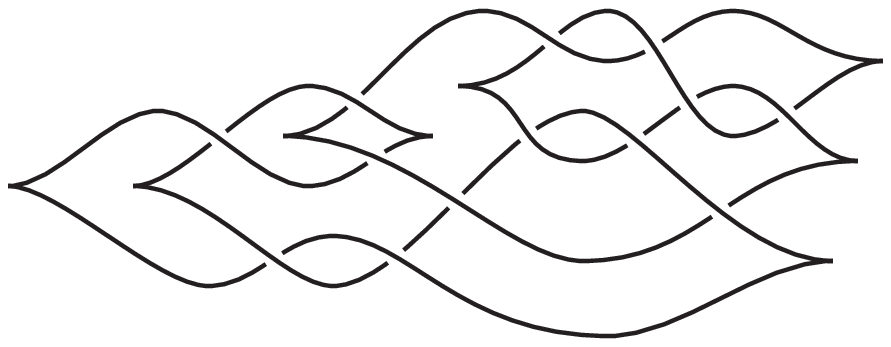}}
\end{Figure}

We consider first decompositions without any 1-handles. In this case, the link
$K_1 \cup \dots \cup K_m \subset \S^3$ can be represented by a
{\sl front projection}, that is a link diagram with horizontal cusps instead of
vertical tangencies, such that at each crossing the arc with most negative slope
crosses in front (cf. figure 1). Then, the Legendrian framing of $K_i$
is given by the {\sl blackboard framing} associated to the diagram with one left-twist
added for each right cusp (see [9]).

In the general case, we represent the 1-handles by dotted circles stacked over
the front projection of a Legendrian tangle, in such a way that the diagram of the link
$K_1, \dots, K_m \subset \#_n\,\S^1 \times \S^2$ is obtained by connecting the
endpoints of the tangle by means of parallel arcs, each one of which pass once through
a dotted circle (cf. figure 2). Again the Legendrian framing of $K_i$ is
given by the blackboard framing associated to the diagram with one left-twist added
for each right cusp.

\begin{Figure}[htb]{}{}{}
\centerline{\fig{Front projection with 1-handles}{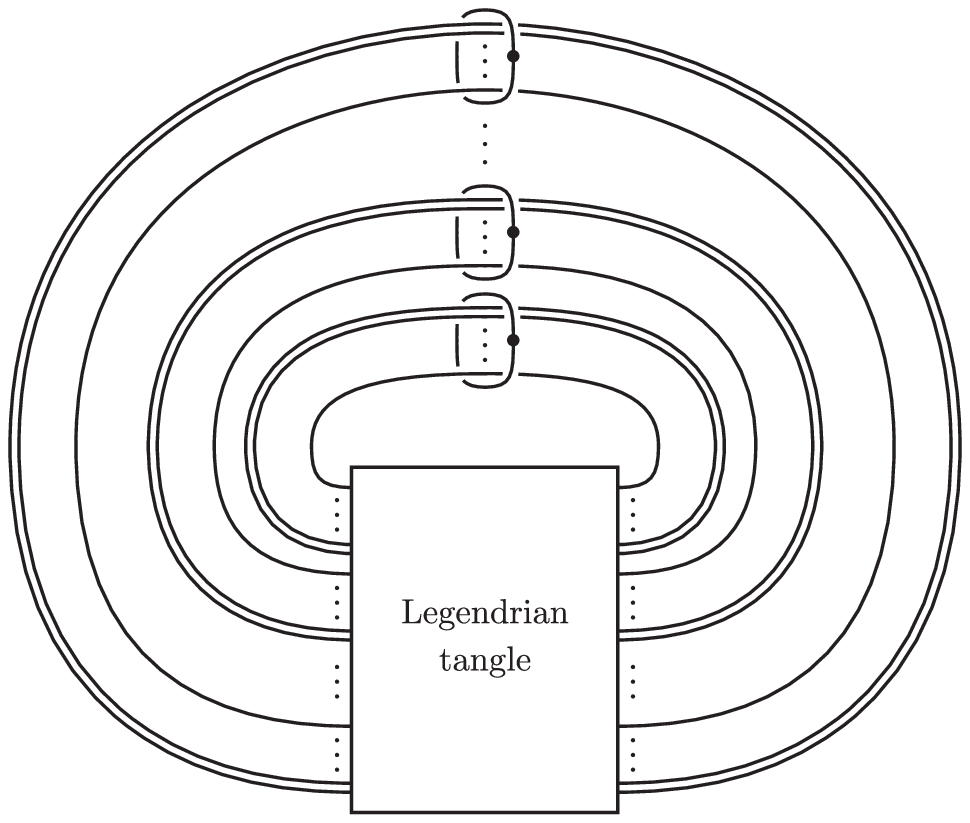}}
\end{Figure}

This way of representing Legendrian 2-handles is the one suggested in \cite{G98},
starting from a {\sl Legendrian link diagram} in {\sl standard form} (cf. definition
2.1 of \cite{G98} and the subsequent discussion at page 634).

In order to get a more convenient representation for our purpose, we modify the handle
decomposition by twisting once negatively each 1-handle. After this change, all the
diagram can be drawn as a front projection with some arcs passing through the dotted
circles, the Legendrian framing still being the blackboard framing with one left-twist
added for each right cusp (cf. figure 3).

\begin{Figure}[htb]{}{}{}
\centerline{\fig{Front projection with twisted 1-handles}{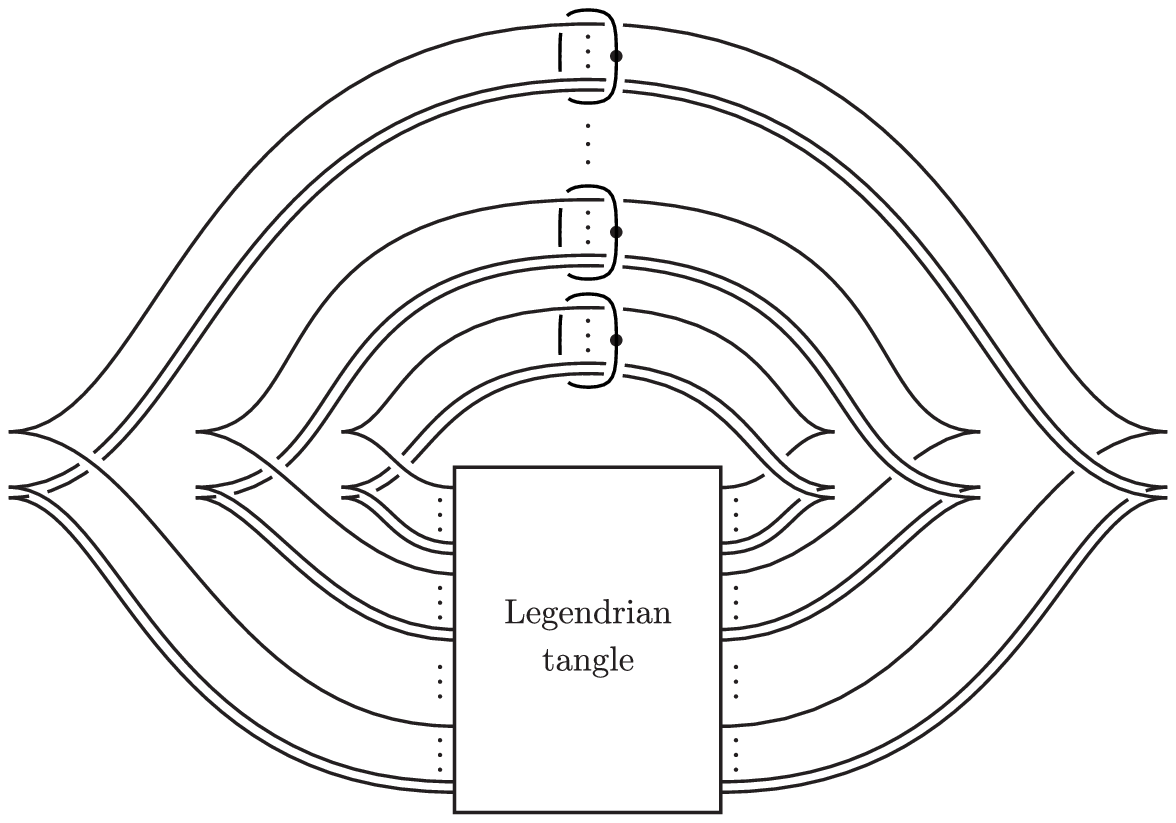}}
\end{Figure}

The following theorem says the all the diagrams considered above do in fact represent
handle decompositions of Stein surfaces. The proof of this fact is implicitely
contained in \cite{E90} (see also \cite{G98}).

\begin{theorem}[Eliashberg] \label{eli/thm}
A smooth oriented compact 4-manifold with boundary is a Stein surface, up to
orientation preserving diffeomorphisms, iff it has a handle decomposition $X_1 \cup
H_1 \cup \dots \cup H_m$, where $X_1$ consists of $\,0$- and $1$-handles and the
$H_i$'s are Legendrian 2-handles attached to $X_1$.
\end{theorem}

Now, we come to the main theorem of this paper, which characterizes compact Stein
surfaces in terms of branched covering and Lefschetz fibration. For proving it,
we will use the fact the any compact Stein surface has a handle decomposition as in
theorem \ref{eli/thm}, but not the viceversa (cf. remark \ref{stein/rem}). 

\begin{theorem} \label{stein/thm}
Given a smooth oriented connected compact $4$-manifold $X$ with boundary,
the following statements are equivalent up to orientation preserving diffeomorphisms:
\begin{enumerate}\itemsep0pt
\item \vskip-\topsep $X$ is a Stein surface;
\item $X$ is an analytic branched covering of $\B^4$;
\item $X$ is a covering of $\B^2 \times \B^2$ branched over a positive braided surface;
\item $X$ is a positive allowable Lefschetz fibration over $\B^2$ with bounded
regular fiber.
\end{enumerate}\vskip0pt\vskip-1\baselineskip\medskip
\end{theorem}

\begin{proof}
(b) $\Rightarrow$ (a). 
Given an analytic branched covering $p:X \to \B^4$, we have that $\Int X$ is a Stein
surface without boundary, since the restriction of $p$ to $\Int X$ is a finite
holomorphic map (see \cite{GR77}, p. 125). Let $f:\Int X \to R$ be a proper strictly
plurisubharmonic function and $g:\Int X \to \R$ be the plurisubharmonic function
defined by $g(x) = 1/(1-\|p(x)\|^2)$. By the transversality of the branch set of $p$
with respect to
$\S^3$, we have $X \cong g^{-1}([0,c])$ for $c > 0$ (regular value) sufficiently large.
Now, the function $h = g + \epsilon f$ is proper and strictly plurisubharmonic on
$\Int X$, for every $\epsilon > 0$. By choosing $\epsilon$ sufficiently small, we have
also $X \cong h^{-1}([0,c])$, hence $X$ is a Stein surface with boundary.

(c) $\Rightarrow$ (b). Let $p:X \to \B^2 \times \B^2$ a covering branched over a
positive  braided surface $S \subset \B^2 \times \B^2$. By theorem \ref{rud/thm}, $p$ is
analitically branched (see \cite{DG94} for the definition). Then, by a theorem of
Grauert and Remmert \cite{GR58} (cf. \cite{DG94}), $p$ is a true analytic covering of
$\B^2 \times \B^2 \cong \B^4$.

(d) $\Rightarrow$ (c). This implication follows immediately from propositions 
\ref{lfB2/thm} and \ref{lfbc/thm}.

(a) $\Rightarrow$ (d). Let $X$ be a Stein surface with boundary. By proposition \ref{bclf/thm},
it is enough to find a simple branched covering $p:X \to \B^2 \times \B^2$, whose
branch set is a positive braided surface. We start with a handle decomposition $X_1
\cup H_1 \cup \dots \cup H_m$,\break where $X_1$ consists of $\,0$- and $1$-handles
and the $H_i$'s are Legendrian 2-handles attached to $X_1$. In order to make the proof 
easier to read, we consider first the special case of one 2-handle attached to $\B^4$.
This allows us to explain the crucial ideas of the proof, avoiding many technical
details. Then, we show how to deal simultaneously with different 2-handles and how to
work the presence of 1-handles.

\smallskip
{\sl Case 1: no 1-handles and one 2-handle.}
In this case, we have $X \cong \B^4 \cup H$, for a Legendrian 2-handle $H$. Let $K
\subset \S^3$ the Legendrian attaching knot of $H$. Then, $K$ can be represented by a
front projection diagram $\cal D$ as described above. An example of such a diagram is
depicted in figure 4; all the diagrams in the following figures 5,
6, 9 and 12 have to be considered as successive modifications of this one. 

\begin{Figure}[htb]{}{}{}
\centerline{\fig{$\cal D$: front projection of a Legendrian kno t}{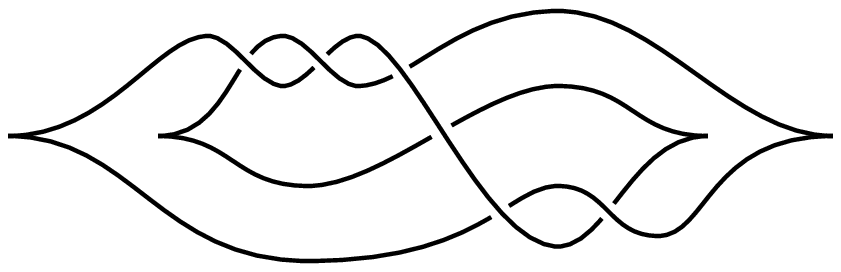}}
\end{Figure}

First of all, we smooth all the cusps and add a negative kink at each right one. In
this way, we get a new diagram $\cal E$ of $K$ (in fact of a transversal knot parallel
to $K$, cf.~\cite{E93}) whose blackboard framing represent the Legendrian framing of
$K$ (see figure 5).

\begin{Figure}[htb]{}{}{}
\centerline{\fig{$\cal E$: negative kinks near to right cusps}{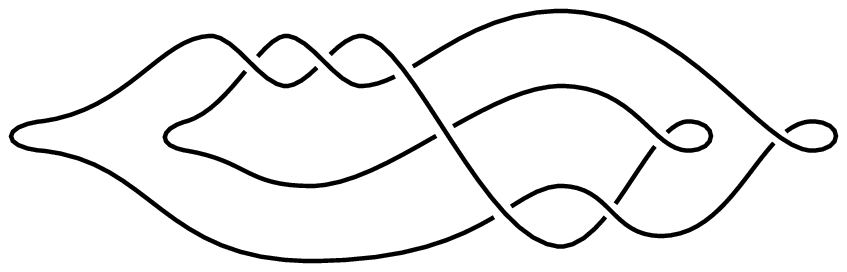}}
\end{Figure}

Then, we redraw $\cal E$ as a polygonal diagram with smoothed corners and edges
of slope $+1$ or $-1$, paying attention to not introduce local minima or maxima for the
abscissa other than the ones coming from cusps, and rotate everything of $-\pi/4$
radians. The resulting diagram $\cal F$ (see figure 6)%
\begin{Figure}[htb]{}{}{}
\centerline{\fig{$\cal F$: horizonal and vertical edges}{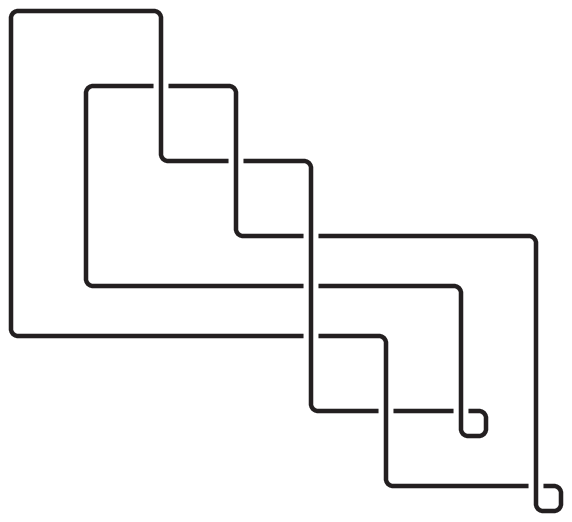}}
\end{Figure} 
has the following properties: all the edges of $\cal F$ are horizontal or vertical; at
each crossing the vertical edge crosses in front; any vertical edge belongs to one of
the three types shown in figure 7, depending on the local structure of
$\cal F$ in a neighborhood of it.

\begin{Figure}[htb]{}{}{}
\centerline{\fig{Vertical edge types 1, 2 and 3}{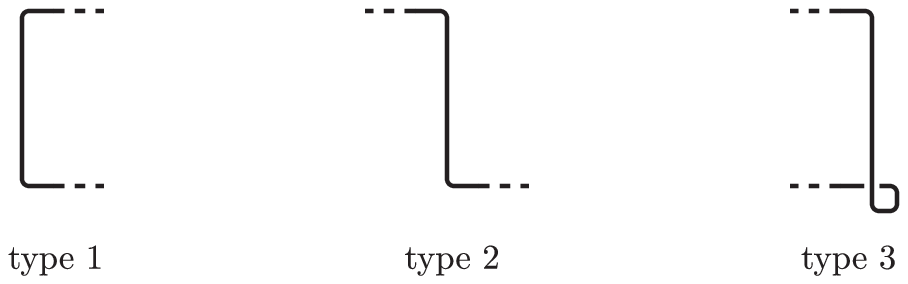}}
\end{Figure}

Finally, we apply to $\cal F$ the moves described in figure 8, in order
to get a new diagram $\cal G$, satisfying the same properties of $\cal F$, with all
the vertical edges of types 1 and 3 respectively in the left-most and the right-most
positions. Of course, also $\cal G$ is a diagram of $K$ (up to smooth
equivalence) whose blackboard framing represent the Legendrian framing of $K$.

\begin{Figure}[htb]{}{}{}
\centerline{\fig{Moving away the vertical edges of type 1 and 3}{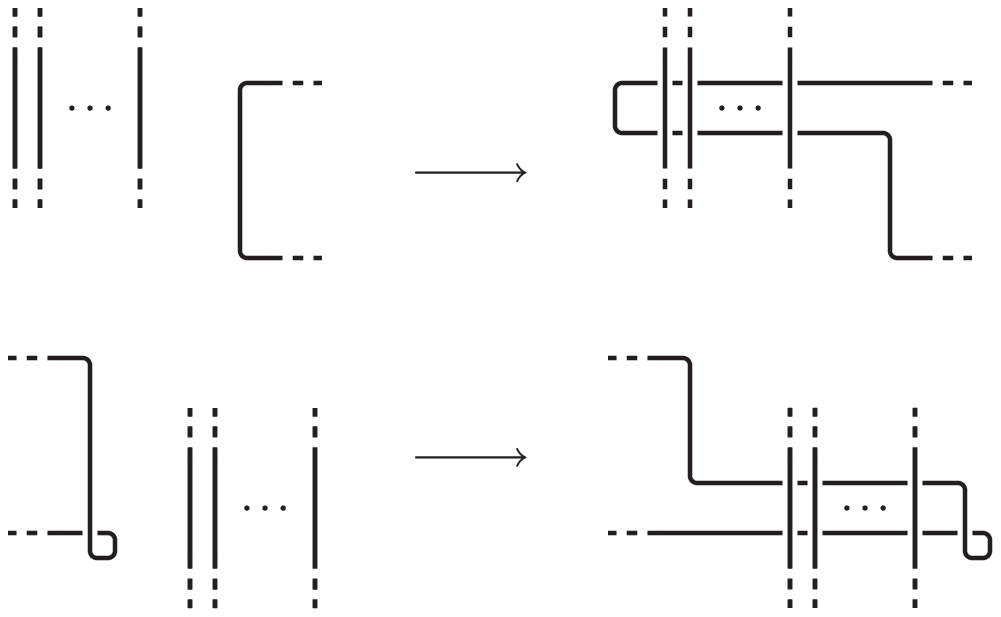}}
\end{Figure}

The vertical edges of the types 1 and 3 come respectively from the
left cusps and the right cusps of the diagram $\cal D$. Hence, putting $c = \#
(\hbox{left cusps of $\cal D$}) = \#(\hbox{right cusps of $\cal D$})$, we have exactly
$c$ vertical edges of type 1 and $c$ vertical edges of type 3.
Let $V_1, \dots, V_{2c}$ be all such edges, numbered starting from the\break
uppermost one of type 1 and following the orientation of the diagram which induces on
it the up-down orientation. We can assume that $\cal G$ has been constructed in such a
way that, going from left to right, we have in the order $V_1,V_3, \dots, V_{2c-1}$ on
the left side of $\cal G$ and $V_2,V_4, \dots, V_{2c}$ on the right side of $\cal G$
(see figure 9).

\begin{Figure}[htb]{}{}{}
\centerline{\fig{$\cal G$: vertical edges of types 1 and 3 moved away}
           {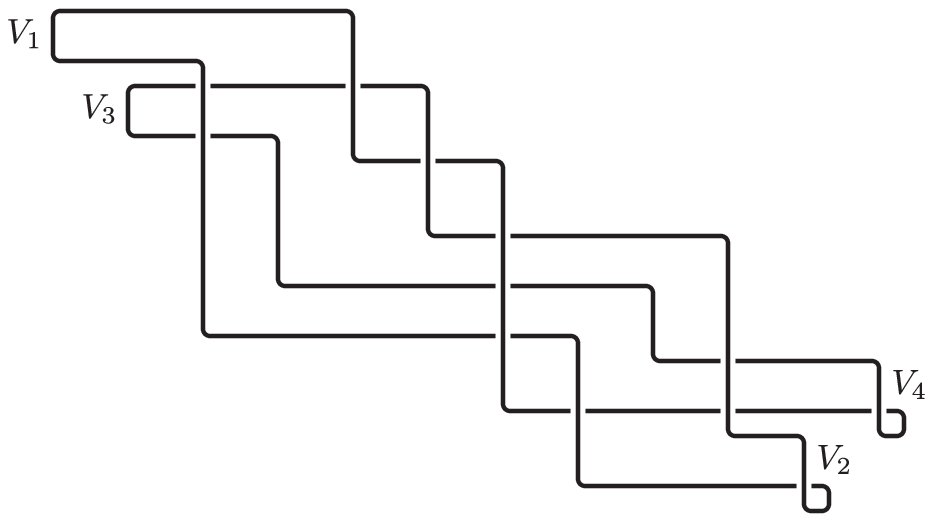}}
\end{Figure}

Now, we consider the simple branched covering $p_0: \B^2 \times \B^2 \to \B^2 \times
\B^2$ with $2c+1$ sheets labelled from $0$ to $2c$, whose branch set consists of
disks $D_1, \dots, D_{2c}$ parallel to the second factor and whose monodromy
around $D_i$ is $(i{-}1\ i)$, for every $i=1, \dots 2c$. We think $D_1, \dots, D_{2c}$
as parallel disks in $\R^3 \subset \S^3 = \Bd \B^4 \cong \B^2 \times \B^2$\break with
interiors pushed inside $\B^4$ and represent their boundaries as vertical lines $L_1,
\dots, L_{2c}$ in the diagram. Furthermore, we assume that: $K \cap D_1 = V_1
\subset L_1$ and $K \cap D_i = \emptyset$ for $i > 1$; $L_i$ lies immediately on the
right (resp. left) of $V_i$ for $i$ odd (resp. even);
$\cal G$ crosses in front of $L_i$ at all the crossings except the upper (resp. lower)
one near to $V_i$ for $i > 1$ odd (resp. even), as shown in figure 10.

\begin{Figure}[htb]{}{}{}
\centerline{\fig{$p_0:\B^2 \times \B^2 \to \B^2 \times \B^2$}{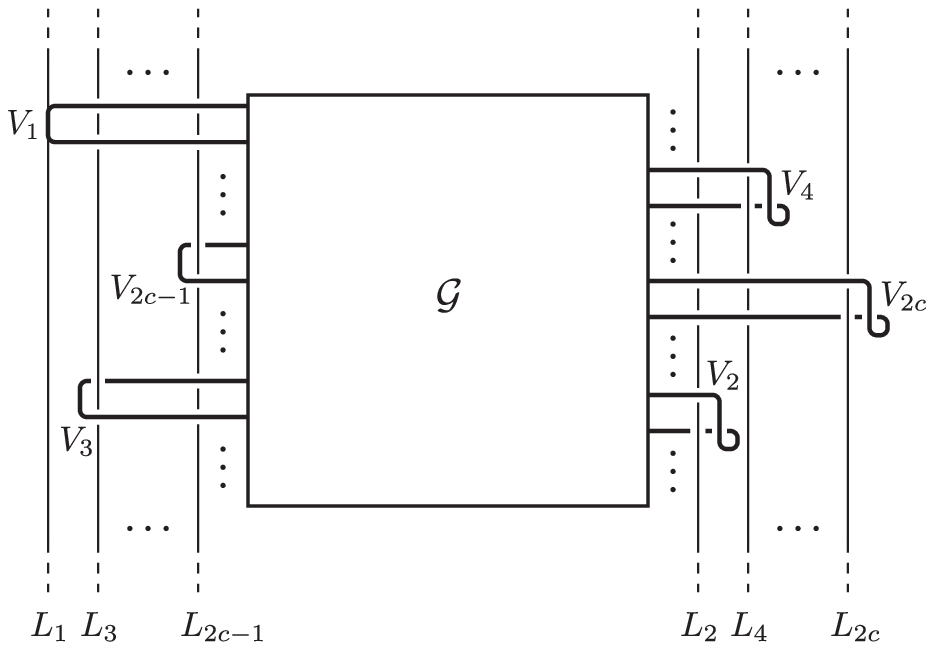}}
\end{Figure}

Let $V_1',\dots, V_{2c}'$ be new vertical edges with the following properties:
$V_i'$ is collinear with $V_i$ for any $i=1,\dots,2c$; all the $V_i'$'s lies above all
the $V_j$'s; the projections of the edges $V_2', \dots ,V_{2c}'$ on $L_1$ have
disjoint interiors and their union coincides with $V_1'$; the bottom end of $V_i'$ and
the top end of $V_{i+1}'$ have the same ordinate for any $i = 2,\dots,{2c-1}$.

Then, we join the $V_i'$'s by horizontal edges, in order to get a trivial knot
diagram linked with the $L_i$'s as shown in figure 11, where the horizontal
edges crosses behind $L_i$ at all the crossings except the lower one near $V_i'$ and
the lowermost one too if $i$ is odd, for any $i > 1$.

\begin{Figure}[htb]{}{}{}
\centerline{\fig{unknot diagram linked to the branch set of $p_0$}{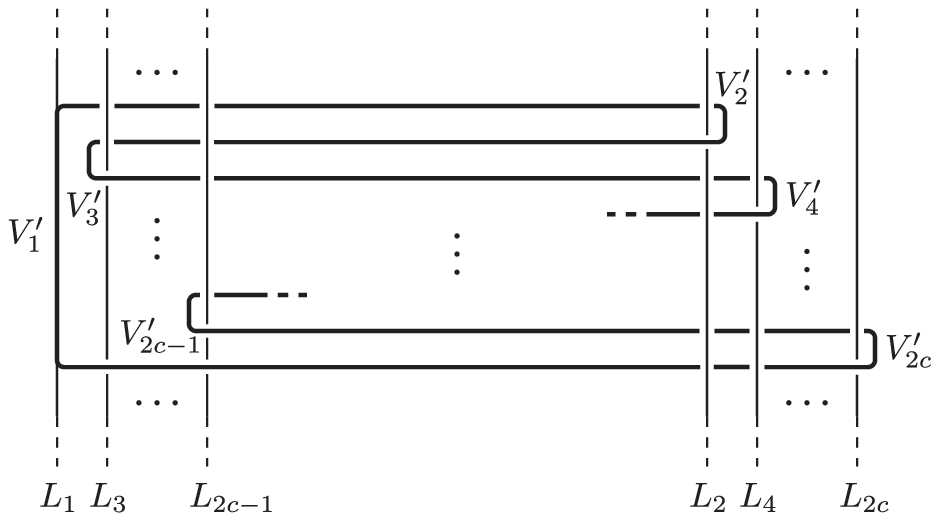}}
\end{Figure}

Finally, we connect this diagram with $\cal G$ by means of a vertical band as show
in figure 12, in such a way that the resulting diagram $\cal H$ is again a
diagram of $K$ intersecting $L_1$ along an arc and the corresponding blackboard
framing still represent the Legendrian framing of $K$.

\begin{Figure}[htb]{}{}{}
\centerline{\fig{$\cal H$ = $\cal G$ $+$ the unknot diagram}{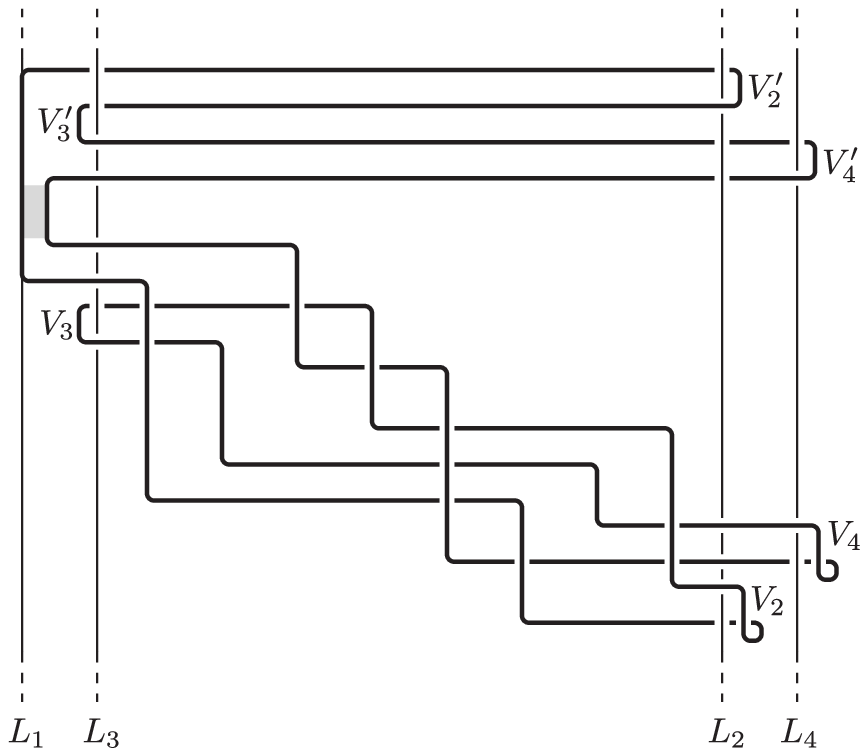}}
\end{Figure}

Let $A \subset K$ be the arc represented by $\Cl(\cal H - L_1)$. Then
$p_0^{-1}(A)$ is the disjoint union of $2c - 1$ arcs and a knot $\widetilde K
\subset \S^3$ equivalent to $K$ by an ambient isotopy of $\S^3$, which makes the
lifting of the blackboard framing along $A$ into the
Legendrian framing of $K$ with one left-twist added. In fact, by unfolding the sheets
of $p_0$ we get a diagram of $\widetilde K$, which is the connected sum of a copy of
$\cal H$ in the  sheet $0$ with a trivial loop going forth and back in the other
sheets. Moreover, the unfolding process, applied to the lifting of the blackboard
framing along $A$, gives us a framing which coincides with the blackboard one except
for a right (resp. left) half-twist for each vertical segment $V_i$ or $V'_i$ with
$i=2,\dots,2c$ odd (resp. even). The knot $\widetilde K$ obtained starting from
figure 12, together with the lifting of the blackboard framing, is represented in
figure 13.

\begin{Figure}[htb]{}{}{}
\centerline{\fig{$\widetilde K$ and the lifting of the blackboard framing}
           {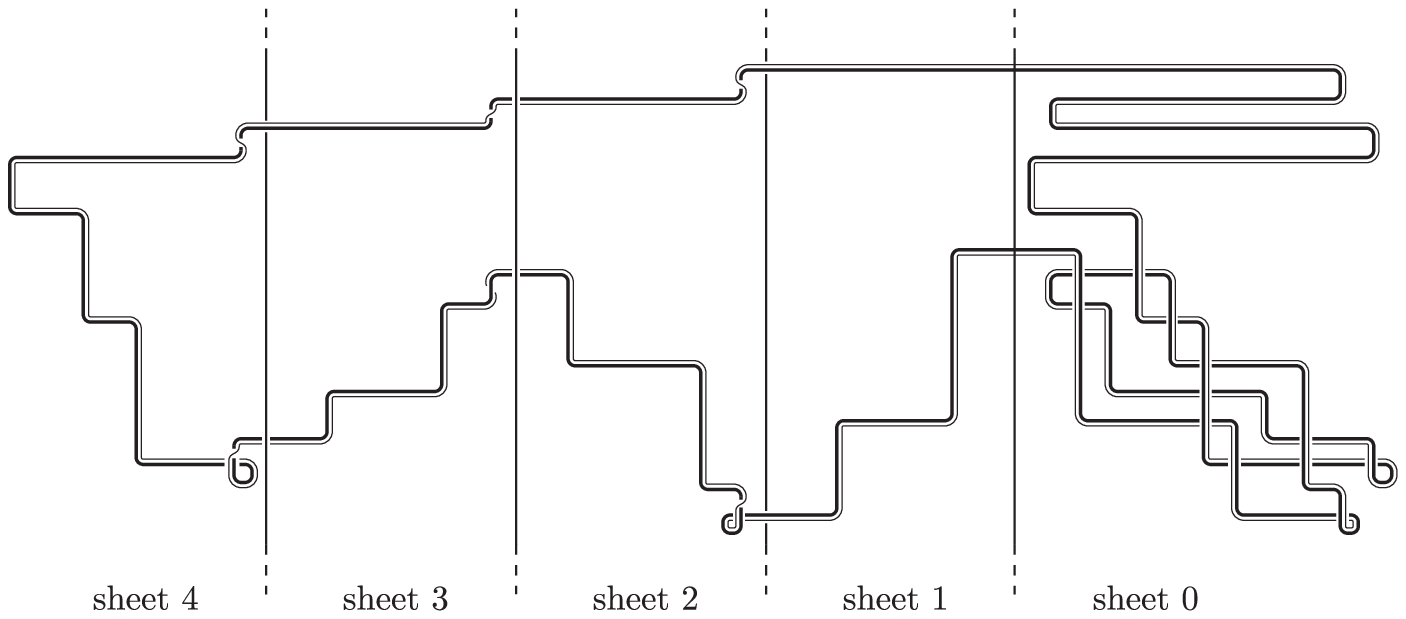}}
\end{Figure}

At this point, the method introduced by Montesinos in \cite{M78} (see also \cite{M80})
allows us to construct a $(2c +1)$-fold simple branched covering $p:M \to \B^2 \times
\B^2$, whose branch set and monodromy coincide with the ones of $p_0$, except for the
attachment to $D_1$ of a ribbon band $B$, which represent the blackboard framing along
$A$ (see figure 14 for the branch set arising from the diagram of
figure 12). Then, denoting by $F_1 \subset \B^2 \times \B^2$ the ribbon
annulus resulting from this surgery on $D_1$, the branch set of $p$ is the regularly
embedded surface
$F_1 \cup D_2 \cup \dots \cup D_{2c} \subset \B^2 \times \B^2$.

\begin{Figure}[htb]{}{}{}
\centerline{\fig{$p:M \to \B^2 \times \B^2$}{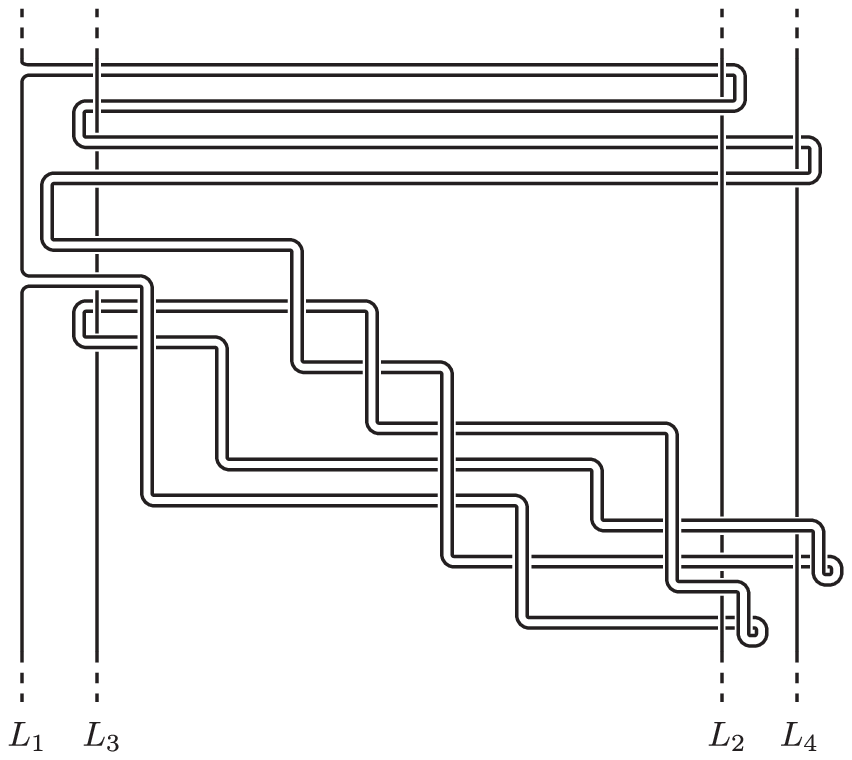}}
\end{Figure}

To conclude this part of the proof, we see that the branch set of $p$ is isotopi\-cally
equivalent to a positive braided surface (over the second factor). 
In fact,\break 
$D_2 \cup \dots \cup D_{2c}$ is already braided (without any twist point) and $F_1$ 
can be made into a braided surface by adapting the Rudolph's braiding process (see
\cite{R83a}) in such a way that all the $D_i$'s are left fixed. Moreover, due to the
special form of $F_1$, all the twist points arising in the process turn out to be
positive.

Namely, we deform the parts of the band $B$ corresponding to vertical edges of $A$ of
types 1, 2 and 3 (including the $V_i'$'s with $i$ odd), one by one from left to right,
to new disks parallel to the $D_i$'s, successively putted in front of the previous
ones, as shown in figure 15.

\begin{Figure}[htb]{}{}{}
\centerline{\fig{The braiding process: vertical edges of types 1, 2 and 3}
           {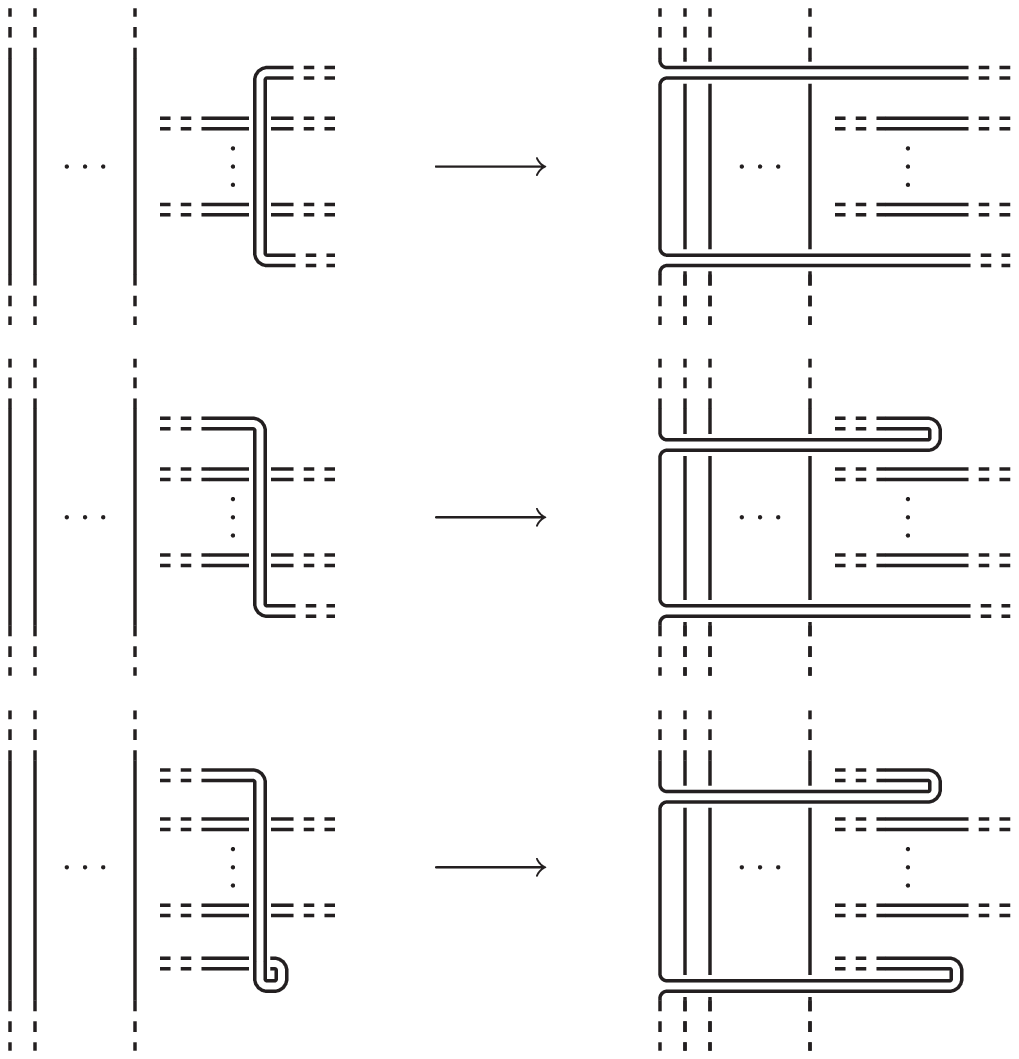}}
\end{Figure}

After all these deformations have been performed, we are left with a certain number of
parallel disks and bands between them (in particular, some of such bands correspond to
the edges $V_i'$ with $i$ even). All such bands have the form depicted in the left part
of figure 16 (up to conjugation), each one being linked to an arbitrary number
(possibly none) of vertical lines. The right part of figure 16, shows how such a band
can isotoped to a braided one with a positive twist point (cf. \cite{R83a}). 

\begin{Figure}[htb]{}{}{}
\centerline{\fig{The braiding process: braiding elementary bands}{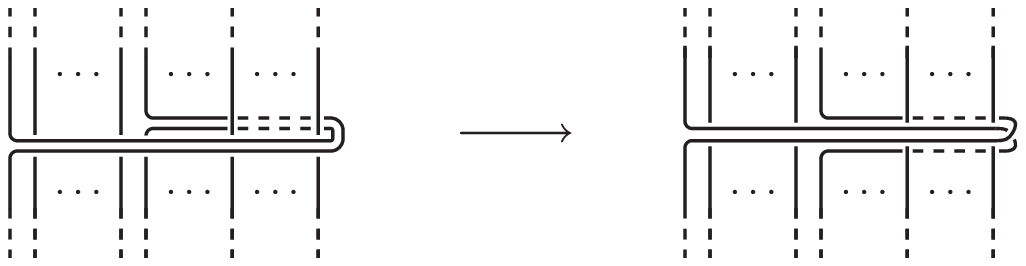}}
\end{Figure}

\smallskip
{\sl Case 2: no 1-handles.}
This time we have $X \cong \B^4 \cup H_1 \cup \dots \cup H_m$, for some Legendrian 
2-handles $H_1, \dots, H_m$. Let $\cal D$ be a front projection diagram of the
Legendrian link $K = K_1 \cup \dots \cup K_m \subset \S^3$, where $K_j$ is the
attanching knot of $H_j$. New diagrams $\cal E$, $\cal F$ and $\cal G$ of $K$ can be
obtained starting from $\cal D$ as in case 1; we use the subscript $j$ for the part of
a diagram corresponding to $K_j$. Then, putting $c_j = \#\hbox{(left cusps of $\cal
D_j$)} = \#\hbox{(right  cusps of $\cal D_j$)}$ and $s_j = c_1 + \dots + c_j$, we
denote by $V_1, \dots, V_{2s_m}$ the vertical edges of types 1 and 3 of $\cal G$.

We assume the $V_i$'s and the $K_j$'s numbered in such a way that: $V_{2s_{j-1}+1},
\dots, V_{2s_j}$ belong to $\cal G_j$ and are ordered in case 1 (starting from the
uppermost of type 1), for any $j=1,\dots,m$; the first edges of the $\cal G_j$'s have
increasing indices from bottom to top, that is we have in the order $V_1, V_{2s_1+1},
\dots, V_{2s_{m-1}+1}$. We also assume the $V_i$'s placed so that, going from left to
right, we have in the order $V_1, V_{2s_1+1}, \dots, V_{2s_{m-1}+1}$\break
$V_3,V_5,\dots,V_{2s_1-1},\dots,V_{2s_1+1}, V_{2s_1+3},\dots,V_{2s_2-1},
\dots,V_{2s_{m-1}+1}$ $V_{2s_{m-1}+3},\dots,V_{2s_m-1}$ on the left side of $\cal G$
and  $V_2,V_4,\dots,V_{2s_m}$ on the right side of $\cal G$.

Then, we consider the simple branched covering $p_0:\B^2 \times \B^2 \to \B^2 \times
\B^2$ with $2s_m+1$ sheets labelled from 0 to $2s_m$, whose branch set consists
of disks $D_1, \dots, D_{2s_m}$ parallel to the second factor and whose monodromy
around $D_i$ is $(0\,\, 2s_j{+}1)$ if $i = 2s_j +1$ and $(i{-}1\ i)$ otherwise.
As above, we think the $D_i$'s as parallel disks in $\R^3$ with the interiors pushed
inside $\B^4$ and represent their boundaries as vertical lines $L_1, \dots, L_{2s_m}$
in the diagram. Furthermore, we assume that: $K \cap D_{2s_{j-1}+1} = V_{2s_{j-1}+1}
\subset L_{2s_{j-1}+1}$ for any $j = 1,\dots, m$; $K \cap D_i = \emptyset$ for all
the other $D_i$'s; the positions of the $L_i$'s and the crossings of $\cal G$ with
them are as in case 1.

Finally, we change each $\cal G_j$ into a new diagram $\cal H_j$, by the same
costruction we have performed in the previous case on the entire diagram $\cal G$ for
obtaining $\cal H$. Thanks to the choices made above about the position of the
$V_i$'s, we can do that without creating any extra crossing. In other words, the new
parts of the diagram, representing the unknots and the bands connecting them with the
$K_j$'s, do not cross each other nor the remaining part of the old diagram $\cal G$. 
Moreover, we let the unknot diagram arising from $\cal G_j$ cross in front of all
the $L_i$'s with $i \neq 2s_{j-1}+1, \dots, 2s_j$.

\begin{Figure}[htb]{}{}{}
\centerline{\fig{$\cal H = \cal H_1 \cup \dots \cup \cal H_m$}{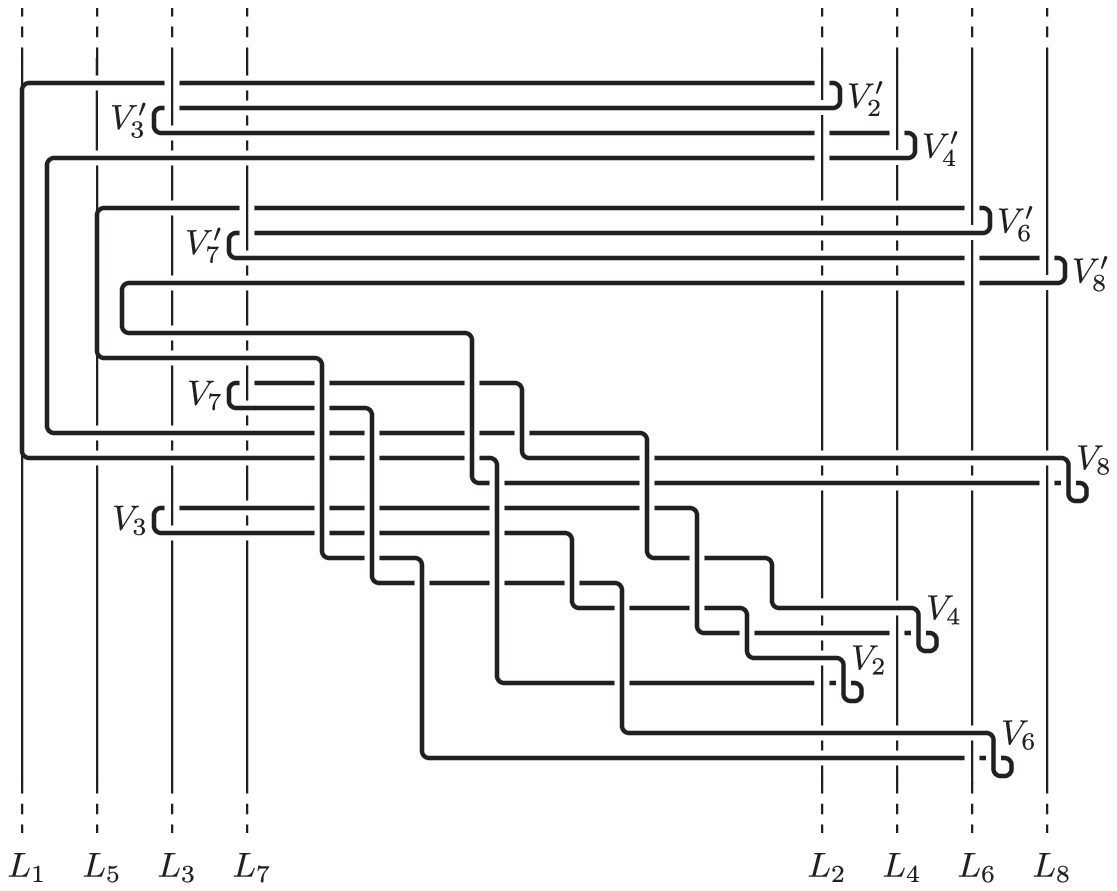}}
\end{Figure}

In this way, we get a new diagram $\cal H = \cal H_1 \cup \dots \cup \cal H_m$ of the
link $K$, such that each $\cal H_j$ meets $L_1 \cup \dots \cup L_{2s_m}$ along an arc
in $L_{2s_{j-1}+1}$ and it is a diagram of $K_j$ whose blackboard framing represent the
Legendrian framing of $K_j$ (see figure 17 for the diagram $\cal H$ obtained
starting with the diagram $\cal D$ of figure 1). 

Let $A = A_1 \cup \dots \cup A_m$, where $A_j \subset K_j$ is the arc represented by
$\Cl(\cal H_j - L_{2s_{j-1}+1})$. Then, $p_0^{-1}(A)$ is the disjoint union of some
arcs and a link $\widetilde K \subset \S^3$ equivalent to $K$ by an ambient isotopy of
$\S^3$, which makes the lifting of the blackboard framing along each $A_j$ into the
Legendrian framing of $K_j$ with one left-twist added. 
We can prove this fact as in case 1, after observing that, as in that case, $\widetilde
K$ is essentially contained in the sheet $0$, being the component $\widetilde K_j$ of
$\widetilde K$ over $K_j$ contained in the sheets $0, 2s_{j-1}+1, \dots, 2s_{j}$, so
that different $\widetilde K_j$'s interact only in the sheet $0$.

In order to get a $(2s_m+1)$-fold simple branched covering $p:M \to \B^2 \times \B^2$,
we modify $p_0$ by attaching to each disk $D_{2s_{j-1}+1}$ a ribbon band $B_j$,
which represets the blackboard framing along $A_j$ and is disjoint from the other
$D_i$'s. Then, the branch set of $p$ is a regularly embedded surface in $\B^2 \times
\B^2$, consisting of $2s_m - m$ disks and $m$ annuli, that can be made into a positive
braided surface, by the same method used in case 1.

\smallskip
{\sl General case.}
Let $X = X_1 \cup H_1 \cup \dots \cup H_m$, where $X_1$ is obtained attaching $n$
1-handles to $\B^4$ and the $H_j$'s are Legendrian 2-handles. We represent such handle
decomposition by a diagram $\cal D$ as in figure 3 and we get diagrams $\cal
E$ and $\cal F$ of $K$ as in the previous cases, expanding the dotted circles behind
the diagram and representing them by dotted vertical lines. So, $\cal F$ crosses in
front of these vertical lines at all the crossings, except the ones corresponding to
passages of the link $K$ through the 1-handles, as shown in figure 18. 

\begin{Figure}[htb]{}{}{} 
\centerline{\fig{$\cal F$ with the dotted vertical lines}{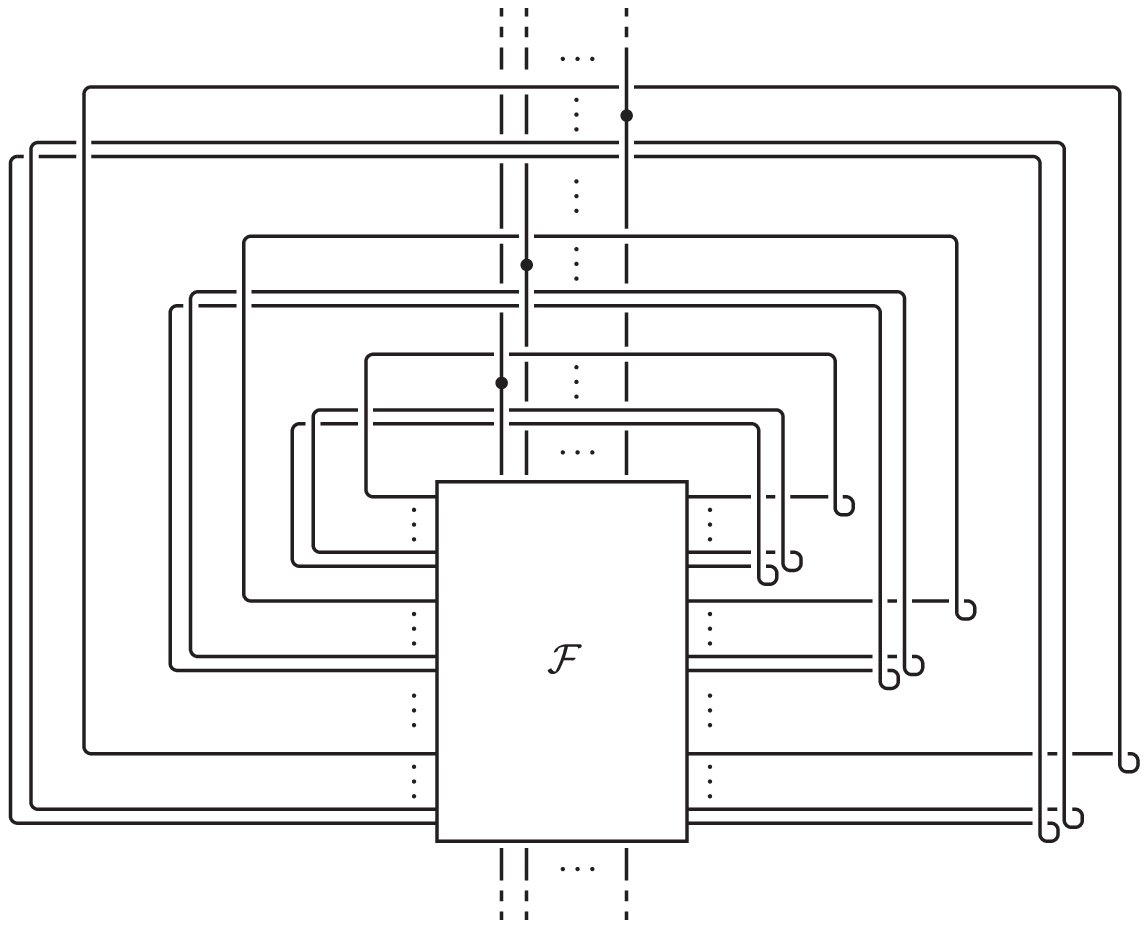}}
\end{Figure}

Then, we push away from $\cal F$ all the vertical edges of type 1 and 3 (including the
ones needed to realize the arcs which go through the 1-handles), by using the moves of
figure 8. In this way, we get a diagram $\cal G$ as in the previous case 2.
We also assume such vertical edges $V_1, \dots, V_{2s_m}$, as well as the subdiagrams
$\cal G_1, \dots, \cal G_m$, numbered and placed as in that case.

Now, let $p_0:\B^2 \times \B^2 \to \B^2 \times \B^2$ the $(2s_m+1)$-simple branched
covering constructed as in case 2, starting from the actual diagram $\cal G$, without
taking into account the dotted components. In order to make $p_0$ into a simple
branched covering $p_1: X_1 \to \B^2 \times \B^2$, we add to it $n$ sheets labelled
from $2s_m +1$ to $2s_m + n$ and $2n$ branch disks $D_{2s_m+1}, \dots, D_{2s_m+2n}$
parallel to the previous ones, whose meridians have monodromies $(0\ 2s_m{+}1),(0\
2s_m{+}1),\dots, (0\ 2s_m{+}\,n), (0\ 2s_m{+}\,n)$.
Assuming also these new disks as parallel disks in $\R^3 \subset \B^2 \times \B^2$
with the interiors pushed inside $\B^4$, we can represent their boundaries in the
diagram by $2n$ vertical lines $L_{2s_m+1}, \dots, L_{2s_m+2n}$. 

We think the $k$-th 1-handle of $X_1$, being realized by the $(2s_m+k)$-th sheet
together with the pair of branch disks $D_{2s_m + 2k-1}$, $D_{2s_m+2k}$ (cf.
\cite{M78}). Then, we draw the lines $L_{2 s_m +2k -1}$ and $L_{s_m +2k}$ in
correspondence of the $k$-th dotted vertical line from the left in figure 18,
letting a horizontal edge of $\cal G$ cross in front of them iff it crosses in front
of such dotted line  (see figure 19).

\begin{Figure}[htb]{}{}{}
\centerline{\fig{$p_1:X_1 \to \B^2 \times \B^2$}{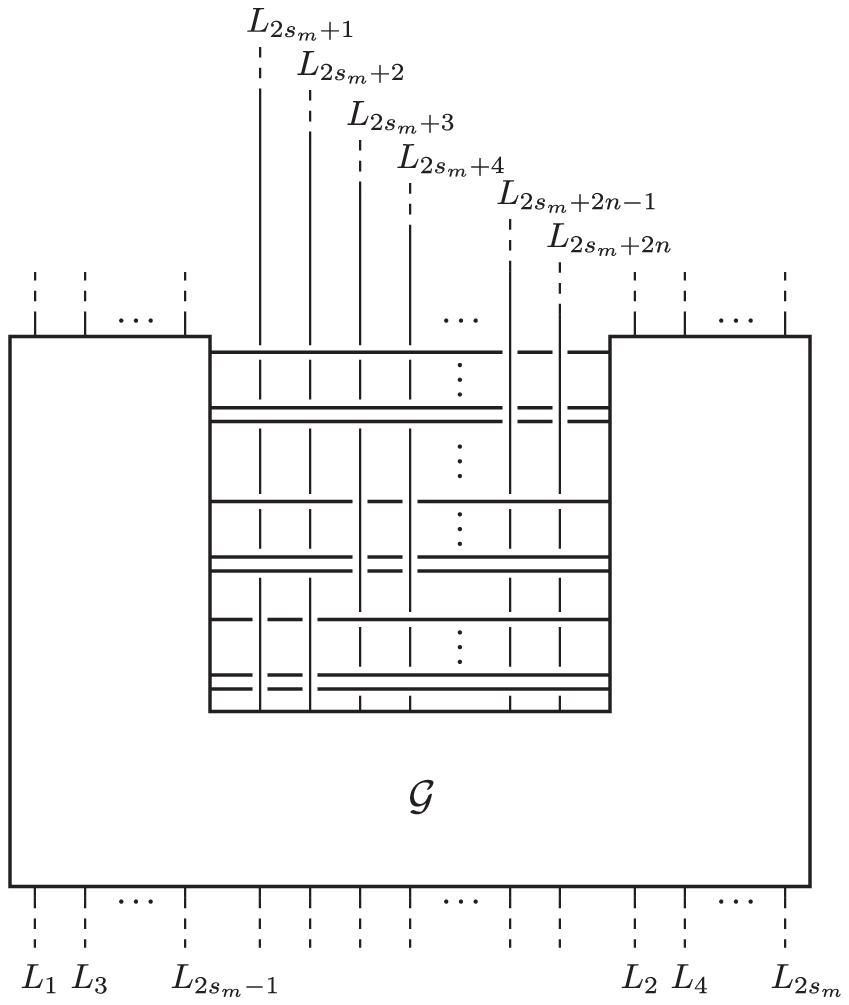}}
\end{Figure}

At this point, we construct another diagram $\cal H$ of $K$, by modifying $\cal G$ as
in case 2 and letting all the new horizontal edges introduced in the construction cross
in front of the vertical lines $L_{2s_m+1}, \dots, L_{2s_m+2n}$.

Finally, we define the disjoint union of arcs $A \subset \S^3$ as above and see, in
the same way, that $p_1^{-1}(A)$ is the disjoint union of some arcs and a link
$\widetilde K \subset X_1$ equivalent to $K$ and that the blackboard framing along each
$A_j$ lifts to the right framing of $\widetilde K_j$. Hence, by attaching to each disk
$D_{2s_{j-1}+1}$ a ribbon band $B_j$ as above, we change $p_1$ into a $(2s_m+n+1)$-fold
simple branched covering $p:M \to \B^2 \times \B^2$. The branch set of $p$ is a
regularly embedded surface in $\B^2 \times \B^2$, consisting of $2s_m + 2n - m$ disks
and $m$ annuli, that can be made into a positive braided surface, again by the same
method used in case 1.
\end{proof} 

\begin{remark} \label{stein/rem}
In proving the implication (a) $\Rightarrow$ (d), we used the hypothesis
only to guarantee the existence of a Legendrian handle decomposition. Then, our proof
of theorem \ref{stein/thm} also provides a new proof of the ``if'' part of
theorem \ref{eli/thm}.

Moreover, we observe that the positivity condition in (c) and (d) is directly related
to the framing properties of Legendrian handles. In fact, by forgetting  such
conditions, we have that: for a 4-manifold as in the statement, having a handle
decomposition with handles of indices $\leq 2$ is equivalent to being  a covering of
$\B^2 \times \B^2$ branched over a braided surface or a Lefschetz fibration over
$\B^2$ with bounded regular fiber (cf. \cite{H79} or \cite{GS99}). 
\end{remark}


\section{Stein fillability\label{fil/sec}}

In this section we apply our main theorem in order to characterize Stein fillable
3-manifolds in terms of open books. First of all, we briefly recall some definitions
and basic facts.

A smooth oriented closed 3-manifold $M$ is called {\sl Stein fillable} iff it is the
oriented boundary of some compact Stein surface $X$ (up to orientation preserving
diffeomorphisms). 
By \cite{B93}, any strictly pseudoconvex boundary of a compact complex surface is Stein
fillable.  Stein fillability is relevant in the context of contact topology of\break
3-manifolds, since the natural contact structure on $M = \Bd X$, given by the complex
tangencies, turns out to be tight (see \cite{E93} or \cite{G98}). The Eliashberg's
characterization of Stein surfaces (theorem \ref{eli/thm}) has been exploited by Gompf in
\cite{G98} for producing several families of fillable 3-manifolds, given in terms of
framed links. 
Using Seiberg-Witten theory, Lisca proved in \cite{L98} that the Poincar\'e homology
sphere with reversed orientation is not Stein fillable (in fact, not simplectically
semi-fillable), as already conjectured in \cite{G98}. Theorem \ref{fil/thm} below, together
with the Harer's  equivalence theorem for fibered links (see \cite{H82}), could enable
us to define an effectively computable obstruction to Stein fillability.

On the other hand, 
given a smooth oriented connected compact surface $F$ with non-empty boundary and a
mapping $\phi \in \Map(F,\Bd F)$, the {\sl open book} with {\sl page} $F$ and {\sl
monodromy} $\phi$ is the space $M_\phi=T(\phi) \cup_k \Bd F$, where $T(\phi)$ is the
mapping torus of $\phi$ and the attaching map $k: T(\phi_{|\Bd F}) \cong \Bd F \times
\S^1 \to \Bd F$ is the projection onto the first factor. It turns out that $M_\phi$ is
a smooth oriented closed 3-manifold (well defined up to orientation preserving
diffeomorphisms) and that $L_\phi = \Bd F \subset M_\phi$ (the {\sl binding} of the
open book) is a fibered link in $M_\phi$ (cf. \cite{H82}). In fact, any such
a 3-manifolds $M$ is orientation preserving diffeomorphic to some open book with
connected binding (see \cite{BE79}).
We say that $M_\phi$ is a {\sl positive} open book iff its monodromy $\phi$ is a
product of right-handed Dehn twists.

The following propositions tell us that the open books coincide, up to orientation
preserving diffeomorphisms, with the boundaries of Lefschetz fibrations over $\B^2$. 

\begin{proposition} \label{lfop/thm}
Let $f:X \to \B^2$ be a  Lefschetz fibration whose regular fiber $F$ has non-empty
boundary. Then $\Bd X$ is orientation preserving diffeomorphic to the open book
$M_\phi$ with page $F$ and monodromy $\phi_f(l) = \phi$, where $l$ is the 
counter\-clockwise loop along $\S^1$.
\end{proposition}

\begin{proof}
Let $y_1,\dots,y_n \in \Int \B^2$ the branch points of $f$ and $l_1,\dots,l_n$
meridian loops around them, such that $l_1 \dots l_n = l$ in $\pi_1(\B^2 - \{y_1,
\dots, y_n\}, \ast)$.
Putting $T = f^{-1}(\S^1)$, we have that the restriction $f_{|T}:T \to \S^1$ is a
locally trivial bundle with fibre $F$ and monodromy $\phi_f\smallcirc i_\ast$,
where $i_\ast$ is the homomorphism indiced by the inclusion of $\S^1$ into the
complement of the branch points $\B^2 - \{y_1,\dots,y_n\}$. Then, $T$ is orientation
preserving diffeomorphic to the mapping torus $T(\phi)$ of the mapping $\phi =
\phi_f(l) = \phi_f(l_1) \dots \phi_f(l_n) \in \Map(F,\Bd F)$.  On the other hand, $T'
= \Cl(\Bd X - T) \cong \B^2 \times \Bd F$, since the restriction $f_{|T'}:T' \to \B^2$
is a (locally) trivial bundle with fiber $\Bd F$. So, we conclude that $\Bd X = T
\cup_{\Bd} T' \cong M_\phi$.
\end{proof}

\begin{proposition} \label{oplf/thm}
For any open book $M_\phi$ with page $F$ there exists a Lefschetz fibration
$f:X \to \B^2$ with regular fiber $F$, such that $\Bd X \cong M_\phi$. Moreover, we can
choose $f$ allowable if $\Bd F$ is connected and positive if $M_\phi$ is
a positive open book.
\end{proposition}

\begin{proof}
Given a open book $M_\phi$ with page $F$, we can write $\phi=
\delta_1^{\epsilon_1}\dots\delta_n^{\epsilon_n}$, with $\delta_i$ right-handed 
Dehn twist along $d_i \subset \Int F$ and $\epsilon _i = \pm 1$. 
Then, fixed $y_1,\dots,y_n \in \Int \B^2$ and $l_1,\dots,l_n$ meridian loops around
them, such that $l_1 \dots l_n = l$ in $\pi_1(\B^2 - \{y_1, \dots, y_n\}, \ast)$, we
consider the Leschetz fibration $f:X \to \B^2$ determined by the branch points $y_1,
\dots,y_n$ and the monodromies $\phi_f(l_i) = \delta_i^{\epsilon_i}$ for $i=1, \dots,
n$\break (cf. section \ref{fib/sec}). By proposition \ref{lfop/thm}, we have 
$\Bd X \cong M$.

For the second part of the proposition, observe that we can choose the $d_i$'s
non-separating if $\Bd F$ connected and the $\epsilon_i$'s positive if
$M_\phi$ is a positive open book. The following lemma \ref{twist/thm} guarantees that
such choices can be made simultaneously. 
\end{proof}

\begin{lemma} \label{twist/thm}
Let $F$ be an oriented connected compact surface with non-empty connected boundary and
let $\delta$ be the right-handed Dehn twist along a simple loop $d \subset \Int F$
parallel to $\Bd F$. Then, there exist right-handed Dehn twists $\delta_1,\dots,
\delta_n$ along non-separating simple loops $d_1,\dots,d_n$, such that $\delta =
\delta_1 \dots \delta_n$ in $\Map(F, \Bd F)$.
\end{lemma}

\begin{proof}
Looking at the double branched covering $p:F \to \B^2$ shown in figure 20, we
see that $d$ covers twice the loop $e \subset \Int \B^2$ encircling all the $2g+1$
branch points, where $g$ denotes the genus of $F$. Then $\delta$ is the lifting of a
double right-handed twist along $e$. By expressing the corresponding braid in terms of
the standard generators, it can be easily realized that $\delta = (\alpha_1\beta_1
\dots \alpha_g\beta_g)^{4g+2}$, where $\alpha_i$ and $\beta_i$ are the right-handed
Dehn twists along the loops $a_i$ and $b_i$ depicted in the figure.
\end{proof}

\begin{Figure}[htb]{}{}{}
\centerline{\fig{$p:F \to \B^2$ and the twist $\delta$}{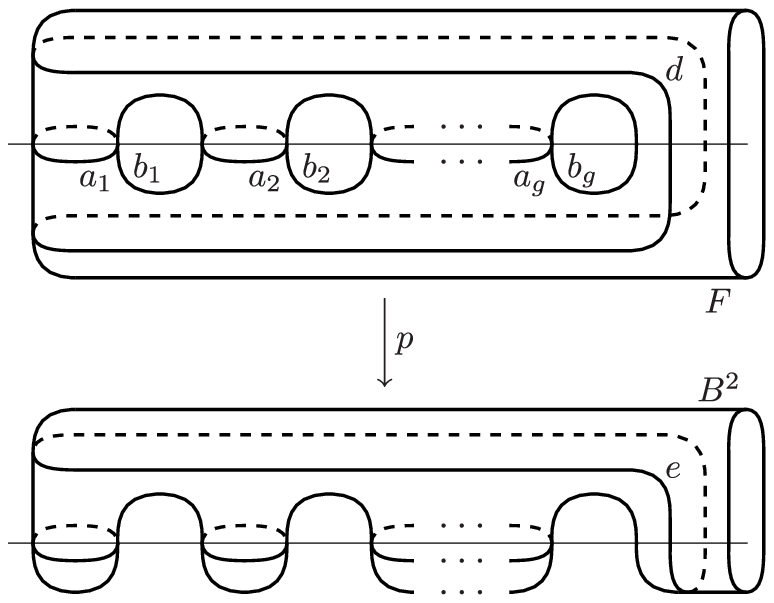}}
\end{Figure}

Now, we are ready to give our fillability criterion.

\begin{theorem} \label{fil/thm}
A smooth oriented closed $3$-manifold is Stein fillable iff it is orientation
preserving diffeomorphic to a positive open-book.
\end{theorem}

\begin{proof}
By theorem \ref{stein/thm} and proposition \ref{lfop/thm}, the oriented boundary 
of any compact Stein surface if orientation preserving diffeomorphic to a positive 
open book. Viceversa, given a positive open book $M_\phi$, we can assume, up to the 
plumbing operation (A) introduced in \cite{H82} (cf. proof of proposition \ref{lfB2/thm} 
above), that the binding of $M_\phi$ is connected. Then, by proposition \ref{oplf/thm} and
theorem \ref{stein/thm}, $M_\phi$ is the oriented boundary of a compact Stein surface. 
\end{proof}

\begin{corollary} \label{surg/thm}
For any smooth oriented closed 3-manifolds $M$ and any fibered knot $K \subset M$,
there is a (possibly trivial) surgery along $K$ which makes $M$ into a Stein fillable
3-manifold.
\end{corollary}

\begin{proof}
Let $M_\phi$ be an open book with page $F$ and binding $L_\phi \subset M_\phi$, such
that $(M,K)$ is orientation preserving diffeomorphic to $(M_\phi,L_\phi)$. Since
$\Map(F,\Bd F)$ is generated by Dehn twists along non-separating simple loops, we can
express $\phi$ as a product of such twists. Now, by lemma \ref{twist/thm}, any
left-handed twist along a non-separating loop can be obtained as a product of some
right-handed twists and of $\delta^{-1}$. In fact, using the notations of
lemma \ref{twist/thm}, this is true for the loop $\delta_1^{-1} =
\delta_2\dots\delta_n\delta^{-1}$, hence the same holds for any non-separating simple
loop in $\Int F$, being all such loops equivalent. Then, we have $\phi = \psi
\delta^{-k}$, with $\psi$ a product of right-handed Dehn twists and $k \geq 0$, because
$\delta$ is a central element of $\Map(F,\Bd F)$. So, we can surgery $M$ along $K$ in
order to get a new 3-manifold $M'$, orientation preserving diffeomorphic to the
positive open book $M_\psi$, which is Stein fillable by theorem \ref{fil/thm}.
\end{proof}


\end{document}